\documentclass[aos,preprint]{imsart}

\usepackage{graphicx,bbm,amssymb,amsmath,amsfonts,amsthm}% ,fullpage,}

\usepackage{stmaryrd}
\usepackage[pdftex]{hyperref}

\def\R{\mathbb{R}}

\def\N{\mathbb{N}}
\def\Z{\mathbb{Z}}
\def\E{\mathbb{E}}
\def\P{\mathbb{P}}

\def\var{\mathrm{Var}}
\def\pen{\mathrm{pen}}
\def\cov{\mathrm{Cov}}
\def\crit{\mathrm{Crit}}

\def\Acal{\mathcal{A}}
\def\Bcal{\mathcal{B}}

\def\Fcal{\mathcal{F}}
\def\Ical{\mathcal{I}}
\def\Mcal{\mathcal{M}}
\def\Tcal{\mathcal{T}}
\def\Xcal{\mathcal{X}}

\def\Ph{\widehat{P}}
\def\ph{\widehat{p}}
\def\tauh{\widehat{\tau}}
\def\muh{\widehat{\mu}}

\def\Pt{\tilde{P}}

\def\Cst{C_{\star}}
\def\Kst{K_{\star}}
\def\Lst{L_{\star}}

\def\taust{\tau_{\star}}
\def\Fcalst{\Fcal_{\star}}

\def\Tbf{\mathbf{T}}

\DeclareMathOperator*{\argmin}{arg\,min}
\DeclareMathOperator{\card}{Card} 

\newcommand{\olP}{\ensuremath{\overline{P}}}

\newtheorem{theorem}{Theorem}
\newtheorem{proposition}[theorem]{Proposition}
\newtheorem{lemma}[theorem]{Lemma}
\newtheorem{definition}[theorem]{Definition}
\newtheorem{corollary}[theorem]{Corollary}

\theoremstyle{remark}
\newtheorem{remark}{Remark}

\newcommand{\paren}[1]{\left( \left. #1 \right. \right)} 
\newcommand{\croch}[1]{\left[ \left. #1 \right. \right]} 
\newcommand{\set}[1]{\left\{ \left. #1 \right. \right\}}
\newcommand{\absj}[1]{\left\lvert #1 \right\rvert} %joli abs
\providecommand{\norm}[1]{\left \lVert #1 \right\rVert}

\newcommand{\PESup}[1]{\left\lceil#1\right\rceil} %Partie Entiere superieure
\newcommand{\PEInf}[1]{\left\lfloor#1\right\rfloor} %Partie Entiere inferieure

\newcommand{\un}{\ensuremath{\mathbf{1}}}

\newcommand{\egaldef}{:=} % egalite definissant la quantite de gauche
\newcommand{\flens}{\rightarrow} % fleche d'application X->Y (ensembles)
\newcommand{\flapp}{\mapsto} % fleche d'application x->f(x) (elements)
\newcommand{\telque}{:}%{\, \mbox{ s.t. } \,} % tel que dans une definition d'ensemble

\newcommand{\hyptag}[1]{\tag{\ensuremath{\mathbf{#1}}}} % nom d'hypothese

\startlocaldefs
\endlocaldefs

\begin{document}

\begin{frontmatter}

% "Title of the paper"
\title{Oracle approach and slope heuristic in context tree estimation}
\runtitle{Oracle approach in context tree estimation}

% indicate corresponding author with \corref{}
% \author{\fnms{John} \snm{Smith}\corref{}\ead[label=e1]{smith@foo.com}\thanksref{t1}}
% \thankstext{t1}{Thanks to somebody} 
% \address{line 1\\ line 2\\ printead{e1}}
% \affiliation{Some University}

\begin{aug}
\author{\fnms{Matthieu} \snm{Lerasle}\thanksref{t1,t2}\ead[label=e1]{mlerasle@unice.fr}} \and
\author{\fnms{Aur\'elien} \snm{Garivier}\thanksref{t2}\ead[label=e2]{aurelien.garivier@telecom-paristech.fr}}

\thankstext{t1}{supported by FAPESP grant 2009/09494-0. This work is part of USP project ``Mathematics, computation, language and the brain''.}
\thankstext{t2}{supported by COFECUB USP 2010 project \emph{Syst\`emes stochastiques \`a interaction de port\'ee variable}}
\runauthor{M. Lerasle and A. Garivier}

\affiliation{IME USP and CNRS Telecom-Paristech}

\address{
Matthieu Lerasle\\
Instituto de Matem\'atica e Estat\'{\i}stica;\\
Universidade de S\~ao Paulo\\
05315-970 S\~ao Paulo, Brasil\\
\printead{e1}\\}

\address{
Aur\'elien Garivier\\
CNRS, Telecom Paristech\\
LTCI UMR 5141; \\
46 rue Barrault 75634 Paris cedex 13\\
\printead{e2}}
\end{aug}

\end{frontmatter}

\begin{abstract}
We introduce a general approach to prove oracle properties in context tree selection. The results derive from a concentration condition that is verified, for example, by mixing processes. Moreover, we show the superiority of the oracle approach from a non-asymptotic point of view in simulations where the classical BIC estimator has nice oracle properties even when it does not recover the source. 

Our second objective is to extend the slope algorithm of \cite{AM08} to context tree estimation. The algorithm gives a practical way to evaluate the leading constant in front of the penalties. We study the slope heuristic underlying this algorithm and obtain the first results on the slope phenomenon in a discrete, non i.i.d framework. We illustrate in simulations the improvement of the oracle properties of BIC estimators by the slope algorithm.
\end{abstract}

\noindent
{\small {{\bf Keywords and phrases :} Context Trees, Penalized Maximum Likelihood, Non-asymptotic Model Selection, Slope heuristic, VLMC, Deviation Inequalities.} }

\section{Introduction}
First motivated by information theoretic considerations, context tree models have been introduced by Rissanen in~\cite{Ri83} as a generalization of discrete Markov models. Since then, they have been widely used in different areas of applied probability and statistics, from  Bioinformatics \cite{bejerano2001a,busch2009} to Linguistics \cite{galves2009,GGGL}.
Sometimes also called Variable Length Markov Chain, a context tree source is a stochastic process whose memory length may vary with the past: the probability distribution of each symbol depends on a finite part of the past, the length of which is a function of the past itself.
Such a relevant part of the past is called a \emph{context}, and the set of all contexts can be represented as a labeled tree called the \emph{context tree} of the process.

Rissanen provided in his seminal paper a pruning algorithm called \emph{Context} for identifying the tree of a process, given a sample $X_1,\dots, X_n$. He proved
his estimator Context to be weakly consistent when the tree of contexts is finite; this result was later completed by a series of papers, including~\cite{BM06} who got rid of the necessity to have a known bound on the maximal length of the memory. On the other hand, penalized maximum likelihood criteria were proved to be strongly consistent in~\cite{CT06,garivier06bicCons}.
More recently, several efforts have been made to obtain non-asymptotic bounds on the probability of correct estimation (see~\cite{Gl11} and references therein).

But the problem of estimation is not the only problem of interest concerning context trees. In fact, these models are widely used because of the remarkable tradeoff they offer between expressivity and simplicity: by providing memory only where necessary, they form a very rich and powerful family of simple processes for the approximation of arbitrary sources. In coding theory, for instance, they are the keystone of the universal coder termed \emph{Context Tree Weighting} (CTW) (see~\cite{willems1995,C01stflour}). The idea behind CTW is that a \emph{double mixture}, over all trees (with a given maximal depth) and, within each tree model, over all parameters, can be computed efficiently. Using this double mixture  as a predictive coding distribution leads to a coder that is proved to satisfy an oracle inequality with respect to the natural loss of Information Theory.
%\[ -\log \nu\left(x_1^n|x_{-\infty}^0\right)\leq \inf_{T}
%\inf_{\theta}-\log
%P_{T,\theta}\left(x_1^n|x_{-\infty}^0\right)+\frac{m-1}{2}|T|\log
%+\frac{n}{|T|} +|T|\left(2+\log
%m \right)+ m-2.\]

The aim of this paper is to show that model selection, and not only aggregation, can be used in an oracle approach for the problem of context tree estimation with K\"ullback loss. 
 For every finite context tree $\tau$ (see Section \ref{sect.notation} for details), we can estimate the transition probabilities of the source $P$ by those $\Ph_{\tau}$ associated to $\tau$ of the empirical measure. The oracle approach consists in looking for the tree minimizing the K\"ullback risk of the estimators $\Ph_{\tau}$. 
 %This approach completes in general the more classical problem of context tree identification, where one wants to recover the minimal tree $\tau_{*}$ of the source. This last problem has been widely studied in the literature, from the initial work \cite{Ri83}, see for example \cite{CT06,GGGL} (mettre des ref). 
This choice of the loss function, while causing a few technical difficulties, emerges naturally from an information theoretic point of view. Following the terminology of~\cite{MN06}, the K\"ullback risk appears as the excess risk associated to the logarithmic loss, which is an (idealized) codelength in coding theory. Hence, the K\"ullback risk appears as a redundancy term caused by the fact that the coder does not know in hindsight which source is to be coded.

When the source has a finite context tree $\tau_{s}$, the oracle approach asymptotically coincides with the consistency approach, because the tree that has the smallest risk is the minimal tree of the source for large numbers of observations. This is no longer the case when the true tree is infinite or at least large compared to the number of data.  Then, the K\"ullback loss of $\Ph_{\tau}$ is decomposed into a bias term measuring the approximation properties of $\tau$ and a variance term measuring the error of estimation. Identification procedures look for the minimal tree with no bias, whereas oracle procedures look for a tree balancing bias and variance. 

The identification approach is inspired from the classical asymptotic situation where the bias term, when non null, is very large compared to the variance term. In this case, under-estimation is easily avoided and the procedures mainly focus on avoiding over-estimation, see for example \cite{CT06, GGGL, Gl11}. On the other hand, the oracle approach is inspired by non asymptotic situations where the true tree is large (compared to the number of data) and can even be infinite. In particular, there exist trees with a bias much smaller than the variance: an oracle is typically a small subtree of $\tau_{s}$ realizing a good tradeoff between the bias (which decreases with the tree size) and the variance (which, in turn,  increases with the tree size). This modern approach is more natural to tackle realistic situations with reasonable number of observations; namely, the set of context trees is used as a toolbox and we want to select the tool that is best suited, in terms of K\"ullback loss. 

The oracle point of view comes from statistical learning theory where it is now well understood in classical problems of non parametric statistics as regression or density estimation (see \cite{Ma07} and the references therein for an introduction). A classical method of selection consists in choosing the model minimizing an empirical loss plus some penalty proportional to the complexity of the model. This principle is the one used in \cite{BBM99, BM97, BM01, Ma07}. Another famous method consists in aggregating a finite set of functions, i.e. to choose a linear combination of previous estimators or approximating functions. An important example of such procedure is the Lasso, where the aggregating weights are chosen by minimization of an $\ell_{1}$-penalized criterion, see \cite{BRT09, CDS01, DET06, GR04, BM06, Ti96, ZH08, ZY07, Zo06}. Complexity penalization procedures are theoretically more interesting because they cover in the same framework several general problems, whereas $\ell_{1}$-penalties are preferred in linear problems for their computational efficiency. We propose a penalization procedure here and we verify that the estimator can be efficiently computed.

Penalized log-likelihood estimators have been studied in context tree estimation, for example by \cite{CT06}. These authors proved that BIC-like estimators (see Section \ref{sect.mod.sel}) are asymptotically consistent when the source has a finite context tree, whatever the leading constant in the BIC-like penalty. Moreover, they showed that BIC estimators can be computed efficiently in practice. However, much less is known about the risk of the selected estimator, when the actual context tree is infinite. In addition, the question of the choice of the leading constant in the BIC penalty for finite number of data remains open. Actually, \cite{GGGL} proved that, for a fixed number of data, the set of trees selected by BIC-like penalties for varying leading constants is exactly the set of champions, where the champion of size $k$ is the tree maximizing the log-likelihood among the trees with less than $k$ degrees of freedom.

Our first goal in this paper is to present a general method to obtain oracle inequalities for a selected $\tauh$, that is, an inequality between the K\"ullback loss of $\Ph_{\tauh}$ and the minimum of the K\"ullback losses of the $\Ph_{\tau}$. We emphasize the central role of concentration inequalities to develop these results for context tree selection, which makes a clear link with model selection theory, as presented, for example, in \cite{BBM99, BM97,BM01} and many others after them. Actually, all the general theorems are consequences of a concentration condition, that we verify for mixing processes. We obtain then a class of examples where the general results apply. For these processes, our penalty takes a BIC form, with a sufficiently large leading constant. As a corollary, we prove therefore that BIC-like estimators have oracle properties when the data are sufficiently mixing.
From a theoretical point of view, the difficulty comes from the fact that new concentration inequalities are required for words that are not contexts, which prevent us from using the martingale approach of~\cite{Cs02,Gl11}.

We study also the slope heuristic of \cite{BM07} in context tree estimation. The heuristic, presented more formally in Section \ref{sect.slope.algo} states the existence of a minimal penalty $\pen_{\min}$ under which the selected tree has huge complexity and over which this complexity is much smaller. Moreover, it states that $2\pen_{\min}$ is an optimal penalty, i.e. that the selected estimator satisfies an asymptotically optimal oracle inequality. The reasons of this phenomenon rely on a fine analysis of the ideal penalty, see \cite{Ar07, AM08}. The ideal penalty is the sum of two terms and the slope heuristic essentially holds when these two terms are asymptotically equal, see \cite{AM08}. The heuristic does not hold in general as was proved in linear regression by \cite{AB10}. In that case, \cite{AB10} proved that an optimal penalty is given by $C\pen_{\min}$ for a constant $C$, different from $2$.   

We study the standard slope heuristic, with an optimal penalty equal to $2\pen_{\min}$, under our concentration assumption and make therefore a contribution to this growing area of statistical learning \cite{AM08, BM07, Le10, Le09}. Note that few proofs are available for non-Hilbertian risks \cite{LT11,Sa11}, and, up to our knowledge, our results are the first ones in a discrete, non i.i.d framework. The heuristic is important since it underlies the slope algorithm presented in \cite{AM08} to calibrate leading constants in the penalties. In the mixing case, the algorithm provides an answer to the question of practical calibration of the leading constant in the BIC penalty. 

We present a simulation study to illustrate our results. 
The simulations are conducted in the particular family of renewal sources (see~\cite{pTHO00a,garivier06redCtw}) for which bias and variance terms can be computed easily, which is not the case in general. 
The simulations show that for relatively small sample sizes of finite sources, the BIC estimator, while failing to recover the true model, does satisfy nice oracle properties; the slope algorithm improves slightly on that, for a very moderately increased computational cost.

The paper is organized as follows. Section \ref{sect.notation} presents some notation used all along the paper. Section \ref{sect.gen.theory} presents our general results. In particular, we show how to deduce from concentration inequalities 1) good penalties yielding oracle properties of the selected estimators and 2) theoretical evidences for the slope heuristic. Section \ref{sect.mix} presents an application of our general approach to mixing processes. We show that they satisfy good concentration properties and we deduce oracle properties of the BIC estimators in this case. Section \ref{sect.simu} presents our simulation study and the proofs are postponed to the appendix.

\section{Notation}\label{sect.notation}
We use the conventions $0/0=+\infty$ and $0\ln(+\infty)=0$. For all $a\in \R$, $\PESup{a}$ denotes the smallest integer larger than or equal to $a$ and $\PEInf{a}$ the largest integer smaller than or equal to $a$. Given two sequences, we use the notation $u_{n}=O(v_{n})$ and $u_{n}=o(v_{n})$ when there exists a constant $C$ such that $|u_{n}|\leq C|v_{n}|$, respectively, when there exists a sequence $\epsilon_{n}\flens 0$ such that $|u_{n}|\leq \epsilon_{n}|v_{n}|$. All the random variables are defined on a probability space $(\Xi,\Xcal,\P)$ and we denote by $\E$ the expectation with respect to $\P$.
Let $A$ be a finite set, with cardinality $|A|$, and, for all $x>0$, let $\log(x)\egaldef \ln(x)/\ln(|A|)$ be the logarithm in base $|A|$. For all $n$ in $\N^{*}$, let $A^{(n)}\egaldef\cup_{k=1,\ldots,n}A^{k}$ and let $A^{*}\egaldef\cup_{k\in \N}A^{k}$. $A$ is called an alphabet and the elements of $A^{*}$ are called words. For all integers $m$ and $n$ such that $m\leq n$, for all words $(a_{m},\ldots,a_{n})\in A^{n-m+1}$, we denote by 
\[
a_{m}^{n}\egaldef(a_{m},\ldots,a_{n})\;\mbox{and}\;\absj{a_{m}^{n}}\egaldef n-m+1\enspace.
\]
The notation $a_{m}^{n}$ is extended to semi-infinite sequences where $m=-\infty$, in that case $a_{-\infty}^{n}\egaldef(a_{i})_{i\leq n}$ and, by definition, for all $n\in \Z$, $\absj{a_{-\infty}^{n}}\egaldef \infty$. The space of semi-infinite sequences is denoted by $A^{-\N}$ and we define $A^{(-\N)}\egaldef A^{-\N}\cup A^{*}$. For every $(\omega,\omega^{\prime})\in A^{(-\N)}\times A^{*}$, let $\omega\omega^{\prime}$ denotes the concatenation of $\omega$ and $\omega^{\prime}$. 
\begin{definition}
A context tree is a subset $\tau\subset A^{(-\N)}$ such that, for every semi-infinite sequence $\omega = a_{-\infty}^{-1}$, there exists a unique $\omega_{\tau}\in \tau$ such that $a_{-\absj{\omega_{\tau}}}^{-1}=\omega_{\tau}$. 

The set of context trees is denoted by $\Tcal$. For every $\tau\in \Tcal$, let $d(\tau)\egaldef\max\set{\absj{\omega},\;\omega\in \tau}$. For every integer $k\geq 1$, let $\Tcal_{k}\egaldef\set{\tau\in \Tcal\telque d(\tau)\leq k}$. When $d(\tau)<\infty$, we say that $\tau$ is finite. For every finite tree $\tau$, let $N(\tau)$ denote the number of elements of $\tau$.

The set $\Tcal$ is provided with the following partial order relation
\[
\tau\prec \overline{\tau}\qquad\mbox{iff}\qquad \forall \omega\in\tau,\;\exists a_{-k}^{-1}\in\overline{\tau}\telque a_{-\absj{\omega}}^{-1}=\omega\enspace .
\]
\end{definition}
In the sequel, we will make repeted use of the following abuse of notation. When $\Acal$ is a set of trees, we will write $\omega\in\Acal$ instead of $\exists \tau\in \Acal \telque \omega\in\tau$. We will in particular use repeatedly the notation $\forall \omega\in\Acal$ instead of $\forall \tau\in \Acal,\;\forall\omega\in\tau$.

\begin{definition}
A transition kernel is a function 
\[P:\Bigg\{\begin{array}{ccc}A^{-\N}\times A & \flens& [0,1]\\(\omega,a)&\flapp& P(a|\omega)\end{array}\] such that, for every $\omega\in A^{-\N}$, $\sum_{a\in A}P(a|\omega)=1$.

A chain $(X_{k})_{k\in \Z}$ is a stationary ergodic stochastic process on $A^{\Z}$.

A chain $(X_{k})_{k\in \Z}$ with distribution $\mu$ on $A^{\Z}$ is said to be \emph{compatible} with transition kernel $P$ if the later is a regular version of the conditional probabilities of the former
\[
\forall (\omega,a)\in A^{-\N}\times A,\qquad \mu\paren{X_{0}=a|X_{-\infty}^{-1}=\omega}=P(a|\omega) \enspace .
\]
\end{definition}

For every chain $(X_{k})_{k\in \Z}$, with distribution $\mu$ compatible with a transition kernel $P$, for every context tree $\tau$, we denote by $P_{\tau}$ a regular version of the following conditional probability:
\[
\forall (\omega,a)\in \tau\times A,\qquad P_{\tau}(a|\omega)\egaldef\mu\paren{X_{0}=a|X_{-\absj{\omega}}^{-1}=\omega}\enspace .
\]
For all finite context trees $\tau$, let $\olP_{\tau}$ be the transition kernel defined by
\[\forall (\omega_{1},\omega,a)\in A^{-\N}\times \tau\times A,\qquad \olP_{\tau}(a|\omega_{1}\omega)\egaldef P_{\tau}(a|\omega)\enspace ,\]
and for all $(\omega,a)\in A^{*}\times a$, let $\mu_{\tau}$ be the probability measure defined recursively by
\[\mu_{\tau}(\omega a)\egaldef\left\{\begin{array}{ll}\mu_{\tau}(\omega) P_{\tau}(a|\omega_{1})&\;\mbox{if}\;\exists (\omega_{1},\omega_{2})\in \tau\times A^{*}\telque \omega=\omega_{2}\omega_{1}\;.\\
\mu(\omega a)&\;\mbox{else.}
\end{array}\right.
\]
\begin{definition}
Let $\Mcal$ be the set of all stationary, ergodic, probability measures on $A^{\Z}$. For all finite context trees $\tau$ let $\Mcal_{\tau}\egaldef\set{\mu\in \Mcal\telque \mu=\mu_{\tau}}$. For every $\mu\in \Mcal_{\tau}$, $(\tau,P_{\tau})$ is called a probabilistic context tree and $\mu$ is called a probabilistic context tree source with tree $\tau$.
\end{definition}
For all transition kernels $Q$, if $\mu\set{ Q(a|\omega)=0\implies P(a|\omega)=0}=1$, we define
\[
K_{\mu}(P,Q)\egaldef \int_{A^{-\N}}d\mu(\omega)\sum_{a\in A}P(a|\omega)\ln\paren{\frac{P(a|\omega)}{Q(a|\omega)}}\enspace .
\]
We take the convention that, if $\mu\set{ Q(a|\omega)=0\implies P(a|\omega)=0}<1$, then $K_{\mu}(P,Q)\egaldef +\infty$. For any finite $\tau$, for any probability measure $\mu$ on $A^{\Z}$ compatible with a transition kernel $P$ and for any family of transition probabilities $(Q(.|\omega))_{\omega\in \tau}$, we also define 
\[
K_{\mu_{\tau}}(P_{\tau},Q)\egaldef \sum_{\omega\in \tau}\mu_{\tau}(\omega)\sum_{a\in A}P_{\tau}(a|\omega)\ln \paren{\frac{P_{\tau}(a|\omega)}{Q(a|\omega)}}\enspace .
\]
The observation set is defined by the projection $X_{1}^{n}$ of a chain $(X_{k})_{k\in\Z}$ with distribution $\mu$ compatible with a transition kernel $P$. Our goal is to estimate $P$ from $X_{1}^{n}$. The risk of the estimators will be measured with the K\"ullback loss $K_{\mu}$. For all $t\leq n$ and all $\omega$ in $A^{(t)}$, we define
\[
\muh_{t}(\omega)\egaldef\frac1{n-\absj{\omega}+1}\sum_{k=\absj{\omega}}^{t}\un_{X_{k-\absj{\omega}+1}^{k}=\omega}\enspace .
\]
A word $\omega$ such that $\muh_{n-1}(\omega)>0$ is called feasible, a tree $\tau$ such that every word is feasible is also called feasible and the set of feasible trees is denoted by $\Fcal$. We also denote, for all $k\leq n$, by $\Fcal_{k}\egaldef\Tcal_{k}\cap \Fcal$.

For all $\tau\in \Fcal$, we denote by $\Pt_{\tau}$ and $\Ph_{\tau}$ the following functions:
\[
\forall (\omega,\omega^{\prime},a)\in \tau\times A^{-\N}\times A,\quad \Ph_{\tau}(a|\omega)=\frac{\muh_{n}(\omega a)}{\muh_{n-1}(\omega)},\quad \Pt_{\tau}(a|\omega^{\prime}\omega)=\Ph_{\tau}(a|\omega)\enspace .
\]
Note that, for all $t\leq n-1$ and $(\omega,a)\in A^{(t)}\times A$, $\sum_{a\in A}\muh_{t+1}(\omega a)=\muh_{t}(\omega)$. Hence, for all feasible trees $\tau$, $\Pt_{\tau}$ defines a transition kernel estimating $P$. Our goal in this paper is to select a tree $\tauh\in\Fcal$ such that, given a confidence level $\delta\in (0,1)$,
\begin{equation}\label{eq.Oracle.1}
\P\set{K_{\mu}\paren{P,\Pt_{\tauh}}\leq \inf_{\tau\in \Fcal^{\prime}} \left[ CK_{\mu}\paren{P,\Pt_{\tau}}+R(\tau,\delta)\right]}\geq 1-\delta\enspace .
\end{equation}
In the previous inequality, the constant $C$ is expected to be close to $1$, the subset $\Fcal^{\prime}\subset \Fcal$ is supposed to be large and the remainder term $R(\tau,\delta)$ should not be too large. In that case, we say that $\tauh$ satisfies an oracle inequality.

Let us mention here that, for every $\tau\in\Fcal$, we have, see Lemma \ref{lem.pyth.gal}, 
\begin{equation}\label{eq.Pythagore}
K_{\mu}\paren{P,\Pt_{\tau}}=K_{\mu}\paren{P,\olP_{\tau}}+K_{\mu_{\tau}}\paren{P_{\tau},\Ph_{\tau}}\enspace .
\end{equation}
In \eqref{eq.Pythagore}, \(K_{\mu}\paren{P,\olP_{\tau}}\) is called the bias term and \(K_{\mu_{\tau}}\paren{P_{\tau},\Ph_{\tau}}\) is called the variance term of the risk. An alternative to \eqref{eq.Oracle.1} is the following
\begin{equation}\label{eq.Oracle.2}
\P\set{K_{\mu}\paren{P,\Pt_{\tauh}}\leq
 \inf_{\tau\in \Fcal^{\prime}} \left[ CK_{\mu}\paren{P,\olP_{\tau}}+R(\tau,\delta)\right]}\geq 1-\delta\enspace .
\end{equation}

%We will first give separate controls of the bias and the variance term of the risk. Then, we will show that BIC estimator satisfy oracle type inequalities \eqref{eq.Oracle.2} when the constant $C$ in front of the penalty is sufficiently large. Finally, we will study the slope heuristic to calibrate this constant in practice.

\section{General approach}\label{sect.gen.theory}

\subsection{Assumptions}
Let us recall the definition of typicality (see \cite{CT06, Cs02} for example).
\begin{definition}
For every $\eta\in(0,1)$ and $k\leq n$, a word $\omega$ is called $(k,\eta)$-typical if 
\[
(1-\eta)\mu(\omega)\leq \muh_{k}(\omega)\leq (1+\eta)\mu(\omega)\enspace .
\]
The set of $(k,\eta)$-typical words is denoted by $\Tbf(k,\eta)$ and let
\begin{equation}\label{def.T.typ}
\Tbf_{typ}(\eta)\egaldef \set{\omega\in A^{(n-1)} \telque \;\forall a\in A,\; \omega\in \Tbf(n-1,\eta)\;\mbox{ and }\;\omega a\in \Tbf(n,\eta)}\enspace .
\end{equation}
\end{definition}

Concentration is the central tool to develop model selection theory, as shown in the series of works \cite{BBM99, BM97,BM01} and many authors after them. In this section, we assume the following concentration condition. There exist $\varrho_{n}\flens 0$, $\varphi_{n}\flens 0$,  $\set{d_{n}\leq n-1, \;d_{n}\flens \infty}$, $\rho_{n}\rightarrow 0$,  and an event $\Omega_{conc}$ satisfying $\P(\Omega_{conc}^{c})\leq \varphi_{n}$, such that
\begin{gather}
\hyptag{CC} \forall \delta\in(0,1),\;\forall t\in\set{n-1,n},\;\forall \tau\in\Fcal_{d_{n}},\;\forall \omega\in\tau\cup\tau\times A\enspace ,\label{hyp.conc}\\
\notag\P\left(\left\{\absj{\muh_{t}(\omega)-\mu(\omega)}\leq \sqrt{\rho_{n}\mu(\omega)\ln\paren{\frac1{\delta}}}+\varrho_{n}\ln\paren{\frac1{\delta}}\right\}\cap\Omega_{conc}\right)\geq 1-\delta\enspace . 
\end{gather}
Let us now choose $d_{n}$ as in assumption \eqref{hyp.conc} and let $\pi$ be a probability measure on $A^{d_{n}+1}$ such that, for all $\omega\in A^{d_{n}}$, $\sum_{a\in A}\pi(\omega a)=\pi(\omega)$. Assumption \eqref{hyp.conc} and a union bound ensure that $\P\set{\Omega_{good}}\geq 1-2\delta-\varphi_{n}$, where
\begin{align}
\notag\Omega_{good}&\egaldef \Bigg\{\forall (\tau, a)\in\Fcal_{d_{n}}\times A, \;\forall\omega\in\tau,\\
&\absj{\muh_{n-1}(\omega)-\mu(\omega)}\leq \sqrt{\rho_{n}\mu(\omega)\ln\paren{\frac1{\pi(\omega)\delta}}}+\varrho_{n}\ln\paren{\frac1{\pi(\omega)\delta}},\nonumber\\
\label{def.omega.good}&\absj{\muh_{n}(\omega a)-\mu(\omega a)}\leq \sqrt{\rho_{n}\mu(\omega a)\ln\paren{\frac1{\pi(\omega a)\delta}}}+\varrho_{n}\ln\paren{\frac1{\pi(\omega a)\delta}}\Bigg\}\; .
\end{align}
Let $\Lambda_{n}^{(1)}\flens \infty$, $\Lambda_{n}^{(2)}\flens \infty$, let $\delta\in (0,1)$, let 
\[
r_{n}(\pi,\delta,\omega,a)=\ln\paren{\frac1{\pi(\omega a)\delta}}\paren{\Lambda_{n}^{(1)}\rho_{n}\vee \Lambda_{n}^{(2)}\varrho_{n}}\enspace ,
\] 
and let 
\begin{align}
\label{def.bon.ens.obs}\Fcalst^{(n)}(\delta)&\egaldef\Big\{ \tau\in \Fcal_{d_{n}}\telque \forall \omega\in\tau, \forall a\in A,\\
&\notag\quad\quad\quad\quad\quad\quad\muh_{n}(\omega a)=0\;\mbox{or}\; \muh_{n}(\omega a)\geq 2r_{n}(\pi,\delta,\omega,a)\Big\}\;.\\
\label{def.bon.ens}\Fcalst(\delta)&\egaldef\Big\{\tau\in \Fcal_{d_{n}}\telque \forall \omega\in\tau, \forall a\in A,\\
&\notag\quad\quad\quad\quad\quad\quad\mu(\omega a)=0\;\mbox{or}\; \mu(\omega a)\geq r_{n}(\pi,\delta,\omega,a)\Big\}\;.\\
\label{def.bon.ens.2}\Fcalst^{(2)}(\delta)&\egaldef\Big\{\tau\in \Fcal_{d_{n}} \telque \forall \omega\in\tau, \forall a\in A,\\
&\notag\quad\quad\quad\quad\quad\quad\mu(\omega a)=0\;\mbox{or}\; \mu(\omega a)\geq 4r_{n}(\pi,\delta,\omega,a)\Big\}\;.
\end{align}
\subsection{A typicality result}
Our first result is that assumption \eqref{hyp.conc} implies typicality of the words in $\Fcalst^{(n)}(\delta)\cup\Fcalst(\delta)$. More precisely, the following proposition holds.
\begin{proposition}\label{prop.typ.gal}
Let $\delta>0$ and let $\Fcalst^{(n)}(\delta)$, $\Fcalst(\delta)$ and $\Fcalst^{(2)}(\delta)$ be the sets defined in \eqref{def.bon.ens.obs}, \eqref{def.bon.ens} and \eqref{def.bon.ens.2} respectively. Let $\Tbf_{typ}(\eta)$ be the set defined in \eqref{def.T.typ} and let $\Omega_{good}$ be the event \eqref{def.omega.good}. There exists $n_{o}$ such that, for all $n\geq n_{o}$, on $\Omega_{good}$, $\Fcalst^{(2)}(\delta)\subset\Fcalst^{(n)}(\delta)\subset\Fcalst(\delta)$. Moreover, there exists  $\eta=O\paren{\sqrt{\frac1{\Lambda_{n}^{(1)}}}\vee \frac1{\Lambda_{n}^{(2)}}}$ such that, on $\Omega_{good}$, all the words in $\Fcalst(\delta)$ belong to $ \Tbf_{typ}(\eta)$.
\end{proposition}
\begin{remark}
Hereafter, we work on the event $\Omega_{good}$ that has ``large probability'', i.e. larger than $1-2\delta-\varphi_{n}$. The collection of trees that we are interested in is $\Fcalst^{(n)}(\delta)$. Proposition \ref{prop.typ.gal} states that, on $\Omega_{good}$, this collection is ``large'' since it contains the collection $\Fcalst^{(2)}(\delta)$ of words with sufficiently large probability of occurrence and the words in $\Fcalst^{(n)}(\delta)$ are typical since they belong to $\Fcalst(\delta)$.
\end{remark}
\subsection{Model selection}\label{sect.mod.sel}
The purpose of this section is to study penalized log-likelihood estimators defined in general as follow. Let $\pen:\Tcal\flens \R_{+}$ and let 
\begin{equation}\label{def.est.pen}
\tauh\egaldef \argmin_{\tau\in \Fcalst^{(n)}(\delta)}\set{\sum_{\omega\in \tau}\muh_{n-1}(\omega)\sum_{a\in A}\Ph_{\tau}(a|\omega)\ln\paren{\frac1{\Ph_{\tau}(a|\omega)}}+\pen(\tau)}\enspace .
\end{equation}
A particular case of such estimators is given by the family of penalties 
\[
\pen_{c}(\tau)=c|A|N(\tau)\frac{\ln n}n\enspace .
\]
This is, up to a constant, the penalty term suggested by the BIC criterion of~\cite{Sc78}: thus, in the following, this penalty will be termed ``BIC-like'' or will even, with some abuse, be called a BIC penalty.
The corresponding estimators have been studied in a series of papers initiated by \cite{CT06} in context tree estimation. It is proved in \cite{CT06} that BIC estimators are consistent when there exists a finite tree $\tau$ such that $\mu=\mu_{\tau}$. We are interested here in oracle properties of the selected estimator, that is, we want to compare $K_{\mu}(P,\Pt_{\tauh})$ with $\inf_{\tau\in \Fcalst^{(n)}(\delta)}K_{\mu}(P,\Pt_{\tau})$. The following theorem is the main result of the paper.
\begin{theorem}\label{theo.oracle.gal}
Let $(X_{n})_{n\in\Z}$ be a stationary ergodic process satisfying assumption \eqref{hyp.conc}. Let $\delta>0$ and let $\Fcalst^{(n)}(\delta)$ be the set defined in \eqref{def.bon.ens.obs}. Let $\Omega_{good}$ be the event \eqref{def.omega.good}. Let 
\[
\ph_{\min}=\min_{(\omega,a)\in\Fcalst^{(n)}(\delta)\times A\telque \Ph(a|\omega)\neq 0}\Ph(a|\omega)\enspace .
\]
Let $L>6+18\ph_{\min}^{-1}$ and let $\tauh$ be the penalized estimator defined in \eqref{def.est.pen}, with
\begin{equation}\label{cond.pen}
\forall \tau\in \Fcalst^{(n)}(\delta),\;\pen(\tau)\geq L\paren{\sqrt{\rho_{n}}+\sqrt{ \frac{\varrho_{n}}{\Lambda_{n}^{(2)}}}}^{2}\sum_{(\omega,a)\in \tau\times A}\ln\paren{\frac1{\pi(\omega a)\delta}}\enspace .
\end{equation}
There exist $n_{o}$ and a constant $\Cst$ such that, for all $n\geq n_{o}$, on $\Omega_{good}$, we have
\[
\forall \tau\in \Fcalst^{(n)}(\delta),\qquad \Cst K_{\mu}(P,\Pt_{\tauh})\leq K_{\mu}(P,\olP_{\tau})+\pen(\tau)\enspace .
\]
\end{theorem}
\begin{remark}
The condition $L>6+18\ph_{\min}^{-1}$ in Theorem~\ref{theo.oracle.gal} can be replaced by  $L>6+36 p_{\min}$, where 
\[p_{\min}=\min_{(\omega,a)\in\Fcalst(\delta)\times A\telque P(a|\omega)\neq 0}P(a|\omega)\enspace ,\]
by using the typicality Proposition~\ref{prop.typ.gal}.
\end{remark}

\begin{remark}
Theorem \ref{theo.oracle.gal} reduces the problem of model selection procedure to the proof of a concentration inequality of the type \eqref{hyp.conc}. We will show in Section \ref{sect.mix} that such concentration inequalities are available when $(X_{n})_{n\in\Z}$ is geometrically $\phi$-mixing. In that case, we can take $d_{n}=O(\ln n)$ and $\rho_{n}=O(n^{-1})$, $\varrho_{n}=O(n^{-1}\ln n)$. Therefore, choosing for $\pi$ the uniform probability measure on $A^{d_{n}+1}$, the condition \eqref{cond.pen} holds for BIC penalties $\pen_{c}(\tau)$ if $c$ is large enough. However, the concentration in that case involves some unknown constant in $\rho_{n}$. Moreover, the constant $L>3+6\ph_{\min}^{-1}$ proposed is too large for practical use. In order to overcome this problem, we propose to study in Section \ref{sect.slope.algo} the slope algorithm of \cite{BM07}.
\end{remark}

\subsection{Slope algorithm}\label{sect.slope.algo}
The slope algorithm has been introduced in \cite{BM07}, it provides a data-driven calibration of the leading constant in a penalty. It is based on the slope algorithm that we adapt here to the particular case of context tree estimation. Let us recall that that the selected tree \(\tauh\) is obtained as a minimizer of the penalized criterion \eqref{def.est.pen}. The heuristic describes the typical behavior of the selected tree \(\tauh\) when \(\pen(\tau)= C\pen_{sh}(\tau)\), $\pen_{sh}(\tau)$ is a well chosen complexity measure of $\tau$ (typically the BIC shape $|A|N(\tau)\frac{\ln n}n$ or the variance term $K_{\mu_{\tau}}(P_{\tau},\Ph_{\tau})$) and $C$ is an increasing leading constant. It states more precisely that there exists a constant $C_{\min}$ such that

\begin{itemize}
\item[{\bf SH1}] When $C<C_{\min}$, the complexity of the selected model $\pen_{sh}(\tauh)$ is very large, typically of the order of $\max_{\tau}\pen_{sh}(\tau)$.
\item[{\bf SH2}] When $C>C_{\min}$, the complexity $\pen_{sh}(\tauh)$ becomes abruptly much smaller.
\item[{\bf SH3}] When $C=2C_{\min}$, the selected estimator satisfies an oracle inequality \eqref{eq.Oracle.1} with a leading constant close to $1$.
\end{itemize}

Let us now assume that we want to calibrate the leading constant $L$ in a penalty of the form $\pen(\tau)=L\pen_{dd}(\tau)$, where $\pen_{dd}(\tau)$ is a data-driven shape for the penalty (typically here, we will use the BIC shape $|A|N(\tau)\frac{\ln n}n$). The slope algorithm evaluates this leading constant in the following data-driven way.

\begin{itemize}
\item[{\bf SA1}] For all $L>0$, compute the complexity $\pen_{dd}(\tauh)$ of the model selected by the penalty $\pen(\tau)=L\pen_{dd}(\tau)$.
\item[{\bf SA2}] Choose $L_{\min}$, such that this complexity is very large for $L<L_{\min}$ and much smaller for $L>L_{\min}$.
\item[{\bf SA3}] Choose finally the constant $L=2L_{\min}$.
\end{itemize}

The algorithm is efficient if, for some constant $L_o$ and some shape penalty $\pen_{sh}$ satisfying the slope heuristic, we have
\[
(L_{o}-o(1))\pen_{dd}(\tau)\leq \pen_{sh}(\tau) \leq (L_{o}+o(1))\pen_{dd}(\tau)\enspace.
\]

Actually, by {\bf SA2}, we observe a jump in the complexity of the selected model when $L\simeq L_{\min}$, hence, by {\bf SH1, SH2},
\[
(L_{\min}-o(1))\pen_{dd}(\tau)\leq C_{\min}\pen_{sh}(\tau) \leq (L_{\min}+o(1))\pen_{dd}(\tau)\enspace.
\]

Therefore, the model selected by {\bf SA3} with the penalty $2L_{\min}\pen_{dd}(\tau)\simeq 2C_{\min}\pen_{sh}(\tau) $ satisfies an oracle inequality, thanks to {\bf SH3}.

The words ``very large'' and ``much smaller'' in Step {\bf SA2}, borrowed from \cite{AM08,BM07}, are not very clear. We refer to \cite{AM08} Section 3.3 for a detailed discussion on what they mean in this context and for precise suggestions on the implementation of the slope algorithm. We refer also to \cite{Michel10} and Section \ref{sect.simu} for practical implementations of the slope algorithm in M-estimation, respectively in our framework.\\

This section presents some theoretical evidences for the slope heuristic. We show the jump of the complexity of the selected model around a minimal penalty as predicted by {\bf SH1}, {\bf SH2} and we prove that a penalty equal to 2 times the minimal one has oracle properties. We do not prove that the leading constant is asymptotically equal to one as predicted by {\bf SH3}. We present finally a theorem that emphasizes what remains to be done to obtain a complete proof of {\bf SH3}. The complexity is $\Delta(\tau)=K_{\mu_{\tau}}(P_{\tau},\Ph_{\tau})$. We prove in Lemma \ref{lem.cons.typ.gal} an upper bound of this term for all $\tau$ in $\Fcalst^{(n)}(\delta)$.

\begin{theorem}\label{theo.min.pen}
Let $(X_{n})_{n\in\Z}$ be a stationary ergodic process satisfying the concentration condition \eqref{hyp.conc}. Let $\delta>0$ and let $\Fcalst^{(n)}(\delta)$ and $\Fcalst(\delta)$ be the sets defined in \eqref{def.bon.ens.obs} and \eqref{def.bon.ens}. Let $\Omega_{good}$ be the event \eqref{def.omega.good}. Let $r>0$ and let $\tauh$ be the penalized estimator defined in \eqref{def.est.pen}, with
\begin{equation}\label{cond.pen.min}
\forall \tau\in \Fcalst^{(n)}(\delta),\qquad 0\leq \pen(\tau)\leq (1-r)K_{\mu_{\tau}}(P_{\tau},\Ph_{\tau})\enspace .
\end{equation}
Let $p_{\min}=\inf_{(\omega,a)\in \Fcalst(\delta)\times A,P(a|\omega)\neq 0}P(a|\omega)$.
Let $\taust$ be the maximizer over $\Fcalst^{(n)}(\delta)$ of $K_{\mu_{\tau}}(P_{\tau},\Ph_{\tau})$ and let $\tau_{o}$ be a minimizer of $K_{\mu}(P,\Ph_{\tau})$ over $\Fcalst^{(n)}(\delta)$. Assume that there exist $\varphi_{MP}\flens 0$ and an event $\Omega_{MP}$ satisfying $\P(\Omega_{MP})\geq 1-\varphi_{MP}$, such that, on $\Omega_{MP}$,
\[
 K_{\mu}(P,\olP_{\taust})=o\paren{ K_{\mu_{\taust}}(P_{\taust},\Ph_{\taust})},\;K_{\mu}(P,\Ph_{\tau_{o}})=o\paren{K_{\mu_{\taust}}(P_{\taust},\Ph_{\taust})}\enspace .
\]
There exists $L\egaldef L(p_{\min},r)$ such that, on $\Omega_{good}\cap\Omega_{MP}$, we have
\begin{equation}\label{eq.min.pen}
K_{\mu_{\tauh}}(P_{\tauh},\Ph_{\tauh})\geq LK_{\mu_{\taust}}(P_{\taust},\Ph_{\taust})\enspace .
\end{equation}
\end{theorem}
\begin{remark}
It is convenient to assume that there exists a constant $p_{o}>0$ such that $p_{\min}\geq p_{o}$. In that case it comes from the proof of Theorem \ref{theo.min.pen} that $L(p_{\min},r)\geq L^{\prime}r$ for some $L^{\prime}$ depending only on $p_{o}$. Theorem \ref{theo.min.pen} states that a penalty smaller than $K_{\mu_{\tau}}(P_{\tau},\Ph_{\tau})$ selects a model with maximal value of $K_{\mu_{\tauh}}(P_{\tauh},\Ph_{\tauh})$. This is exactly {\bf SH1} for the complexity measure $\Delta(\tau)=K_{\mu_{\tau}}(P_{\tau},\Ph_{\tau})$. 
\end{remark}
\begin{remark}
The extra assumption $K_{\mu}(P,\olP_{\taust})=o\paren{ K_{\mu_{\taust}}(P_{\taust},\Ph_{\taust})}$ is natural, since the model with maximal complexity is likely to have a lot of leaves, and therefore a small bias. 
\end{remark}
\begin{remark}
The assumption $K_{\mu}(P,\Ph_{\tau_{o}})=o\paren{K_{\mu_{\taust}}(P_{\taust},\Ph_{\taust})}$ means that the risk of an oracle is much smaller than the maximal risk. It is a natural assumption for the slope heuristic to hold, actually, in {\bf SH2, SH3}, the interesting models are those with a complexity much smaller than the biggest one.
\end{remark}

\begin{theorem}\label{theo.opt.pen}
Let $(X_{n})_{n\in\Z}$ be a stationary ergodic process satisfying the concentration condition \eqref{hyp.conc}. Let $\delta>0$ and let $\Fcalst^{(n)}(\delta)$ be the set defined in \eqref{def.bon.ens.obs}. Let $\Omega_{good}$ be the event \eqref{def.omega.good}. Let $r_{1}>0$, $r_{2}>0$ and let $\tauh$ be the penalized estimator defined in \eqref{def.est.pen}, with
\begin{equation}\label{cond.pen.good}
\forall \tau\in \Fcalst^{(n)}(\delta),\qquad (1+r_{1})K_{\mu_{\tau}}(P_{\tau},\Ph_{\tau})\leq \pen(\tau)\leq (1+r_{2})K_{\mu_{\tau}}(P_{\tau},\Ph_{\tau})\enspace .
\end{equation}
Let $\taust$ be the maximizer of $K_{\mu_{\tau}}(P_{\tau},\Ph_{\tau})$ over $\Fcalst^{(n)}(\delta)$ and let $\tau_{o}$ be a minimizer of $K_{\mu}(P,\Pt_{\tau})$ over $\Fcalst^{(n)}(\delta)$. 
 Assume that there exist $\varphi_{MP}\flens 0$ and an event $\Omega_{MP}$ satisfying $\P(\Omega_{MP})\geq 1-\varphi_{MP}$, such that, on $\Omega_{MP}$,
\[
 K_{\mu}(P,\olP_{\taust})=o\paren{ K_{\mu_{\taust}}(P_{\taust},\Ph_{\taust})},\;K_{\mu}(P,\Ph_{\tau_{o}})=o\paren{K_{\mu_{\taust}}(P_{\taust},\Ph_{\taust})}\enspace ,
\]
\[
\forall \tau\in \Fcalst^{(n)}(\delta),\qquad \sqrt{K_{\mu}(P,\olP_{\tau})K_{\mu_{\tau}}(P_{\tau},\Ph_{\tau})}=o\paren{K_{\mu_{\taust}}(P_{\taust},\Ph_{\taust})}\enspace .
\]
On $\Omega_{good}\cap\Omega_{MP}$, we have
\begin{equation}\label{eq.opt.pen.res.1}
K_{\mu_{\tauh}}(P_{\tauh},\Ph_{\tauh})=o\paren{K_{\mu_{\taust}}(P_{\taust},\Ph_{\taust})}\enspace .
\end{equation}
In addition, a sequence $\eta=O\paren{\sqrt{\frac1{\Lambda_{n}^{(1)}}}\vee \frac1{\Lambda_{n}^{(2)}}}$ exists such that, $\forall \epsilon_{1}>0$, $\forall \epsilon_{2}>0$, there exists $\Cst\egaldef \Cst(r_{1},r_{2},\epsilon_{1},\epsilon_{2},p_{\min})$ such that, on $\Omega_{good}\cap\Omega_{MP}$,
\begin{equation}\label{eq.opt.pen}
\croch{(1-\epsilon_{1}-\epsilon_{2})\wedge \paren{r_{1}-\epsilon_{1}-\frac{\Lst^{(\tau_{o},\tauh)}}{\epsilon_{2}}-\eta}}K_{\mu}(P,\Pt_{\tauh})\leq \Cst K_{\mu}(P,\Ph_{\tau_{o}})\enspace .
\end{equation}
\end{theorem}

\begin{remark}
The assumption $\forall \tau\in \Fcalst^{(n)}(\delta),\; \sqrt{K_{\mu}(P,\olP_{\tau})K_{\mu_{\tau}}(P_{\tau},\Ph_{\tau})}=o\paren{K_{\mu_{\taust}}(P_{\taust},\Ph_{\taust})}$ means that there is no model with a lot of bias and a big variance. It typically holds when trees with large variance are those with a lot of leaves whereas trees with a large bias are the small ones.
\end{remark}
\begin{remark}
{\bf SH2} immediately follows from \eqref{eq.opt.pen.res.1} since, as soon as the penalty becomes larger than $K_{\mu_{\tau}}(P_{\tau},\Ph_{\tau})$ the complexity of $\tauh$ becomes much smaller than the largest one. {\bf SH3} follows partially from \eqref{eq.opt.pen}. The oracle property implies in particular the convergence of the quantities $P_{\tauh}(a|\omega)$ to $P_{\tau_{o}}(a|\omega)$ so that $\Lst^{(\tau_{o},\tauh)}\rightarrow 1$. Therefore, the condition $r_{1}>\Lst^{(\tau_{o},\tauh)}$ to obtain the oracle inequality became asymptotically $r_{1}>1$, and the condition \eqref{cond.pen.good} on the penalty becomes $\pen(\tau)>2\pen_{\min}(\tau)$. \eqref{eq.opt.pen} states then that $2\pen_{\min}$ is asymptotically a penalty yielding an oracle inequality, but this is not exactly {\bf SH3} which states, moreover, that $\Cst \rightarrow 1$. 
\end{remark}

As mentioned, we did not completely prove point {\bf SH3} of the heuristic. The following theorem emphasizes the missing point of the proof. In order to state the result, let us define, for all $\tau\in \Fcal$,
\[
L(P_{\tau})\egaldef \sum_{(\omega,a)\in \tau\times A}(\muh_{n}(\omega a)-\mu(\omega a))\ln \paren{\frac{1}{P_{\tau}(a|\omega)}}\enspace .
\]

\begin{theorem}\label{theo.opt.pen.afaire}
Let $(X_{n})_{n\in\Z}$ be a stationary ergodic process satisfying the concentration condition \eqref{hyp.conc}. Let $\delta>0$ and let $\Fcalst^{(n)}(\delta)$ be the set defined in \eqref{def.bon.ens.obs}. Let $\Omega_{good}$ be the event \eqref{def.omega.good}. Let $r_{1}>0$, $r_{2}>0$ and let $\tauh$ be the penalized estimator defined in \eqref{def.est.pen}, with a penalty term satisfying \eqref{cond.pen.good}.
Let $\tau_{o}$ be a minimizer of $K_{\mu}(P,\Pt_{\tau})$ over $\Fcalst^{(n)}(\delta)$. 
 Let $u<1,v<1$ and 
\[
\Omega_{mis}\egaldef \set{\forall \tau\in \Fcalst^{(n)}(\delta),\;L(P_{\tau})-L(P_{\tau_{o}})\leq uK_{\mu}(P,\Ph_{\tau})+vK_{\mu}(P,\Ph_{\tau_{o}})}\enspace .
\]
On $\Omega_{good}\cap\Omega_{mis}$, there exists $\eta=O\paren{\sqrt{\frac1{\Lambda^{(n)}_{1}}}\vee \frac1{\Lambda^{(n)}_{2}}}$ such that,
\begin{equation}\label{eq.tout.va.bien}
(1-u)K_{\mu}(P,\olP_{\tau})+(1+r_{1}-u-\eta)K_{\mu_{\tauh}}(P_{\tauh},\Ph_{\tauh})\leq (1+r_{2}+v+\eta)K_{\mu}(P,\Ph_{\tau})\enspace .
\end{equation}
\end{theorem}
\begin{remark}
Assume that there exist $u\flens 0$, $v\flens 0$, $\varphi\flens 0$ such that $\P\set{\Omega_{mis}}\geq 1-\varphi$. Then, from \eqref{eq.tout.va.bien}, for any $r_{1}>0$, the complexity of the selected model is the one of an oracle, that should be much smaller than the maximal one as already explained. This is {\bf SH2} with $\pen_{\min}(\tau)=K_{\mu_{\tau}}(P_{\tau},\Ph_{\tau})$. Moreover, for $r_{1}=r_{2}=1$, i.e. $\pen(\tau)=2K_{\mu_{\tau}}(P_{\tau},\Ph_{\tau})=2\pen_{\min}(\tau)$, \eqref{eq.tout.va.bien} shows that the risk of the selected model is asymptotically exactly the one of an oracle. This is {\bf SH3}.
\end{remark}
\begin{remark}
The weakness of Theorem \ref{theo.opt.pen} comes from the fact that we were not able with our approach to prove that $\Omega_{mis}$ holds with large probability with $u, v\flens 0$. We only obtain this result for some $u>0$, $v>0$. 
\end{remark}
\begin{remark}
In the mixing case that we develop in Section \ref{sect.mix}, we can show that $\Omega_{mis}$ holds with large probability with $u, v\flens 0$ if there exists a fixed set $\Fcalst$ containing $\Fcalst^{(n)}(\delta)$ with large probability such that $\card\set{\Fcalst}=O(n^{\alpha})$ see Proposition \ref{prop.slope.mix} in Section \ref{section.slope.mix}.
\end{remark}
\begin{remark}
Theorem \ref{theo.opt.pen.afaire} shows a difference between context tree estimation and other classical problems of regression or density estimation, where the slope heuristic has been proved. In these frameworks, it is easy to prove that $\Omega_{mis}$ holds with large probability as a consequence of Benett's concentration inequality, see \cite{AM08, Le09}. The main difficulty for proving the slope heuristic is then to show that, with our notation, $K_{\mu_{\tau}}(P_{\tau},\Ph_{\tau})\simeq K_{\muh}(\Ph_{\tau},P_{\tau})$ see for example \cite{AM08}. In context tree estimation, this last result is a direct consequence of typicality, as shown by Lemma \ref{lem.cons.typ.gal} and the problem of $\Omega_{mis}$ seems harder. 
\end{remark}

\section{Application in the mixing case}\label{sect.mix}

We showed in the previous section that oracle inequalities and the slope heuristic can be derived from the concentration condition \eqref{hyp.conc}. Our aim in this section is to show that such concentration result holds for mixing processes.

Let us recall the definition of $\beta$-mixing and $\phi$-mixing coefficients, due respectively to \cite{RV59} and \cite{Ib62}. Let $(\Xi,\Xcal,\P)$ be a probability space and let $\Acal$ and $\Bcal$ be two $\sigma$-algebras included in $\Xcal$. We define
\begin{align*}
\beta(\Acal,\Bcal)&=\frac12\sup\set{\sum_{i=1}^{I}\sum_{j=1}^{J}\absj{\P\set{A_{i}\cap B_{j}}-\P\set{A_{i}}\P\set{B_{j}}}}\enspace,\\
\phi(\Acal,\Bcal)&=\sup_{A\in\Acal,\;\P\set{A}>0}\sup_{B\in\Bcal}\set{\P\set{B|A}-\P\set{B}}\enspace .
\end{align*}
The first $\sup$ is taken among all the finite partitions of $\Xi$ $(A_{i})_{i=1,\ldots,I}$ and $(B_{j})_{j=1,\ldots,J}$  such that, for all $i=1,\ldots,I$, $A_{i}\in \Acal$ and for all $j=1,\ldots,J$, $B_{j}\in \Bcal$.\\
For all stationary sequences of variables $(X_n)_{n\in\Z}$ defined on $(\Xi,\Xcal,\P)$, let
\begin{equation*}
\beta_k=\beta(\sigma(X_i,  i\leq 0),\sigma(X_i,  i\geq k)),\qquad \phi_{k}=\phi(\sigma(X_i,  i\leq 0),\sigma(X_i,  i\geq k))\enspace .
\end{equation*}
The process $(X_n)_{n\in\Z}$ is said to be $\beta$-mixing when $\beta_k\rightarrow 0$ as $k\rightarrow \infty$, it is said to be $\phi$-mixing when $\phi_k\rightarrow 0$ as $k\rightarrow \infty$. It is easy to check, see for example inequality (1.11) in \cite{Br05}, that $\beta(\Acal,\Bcal)\leq \phi(\Acal,\Bcal)$ so that $(X_n)_{n\in\Z}$ is $\phi$-mixing implies $(X_n)_{n\in\Z}$ is $\beta$-mixing.

\begin{theorem}\label{theo.conc.bas}
Let $(X_n)_{n\in\Z}$ be a $\phi$-mixing process satisfying
\begin{equation}\label{cond.phi.mix}
\hyptag{MC}\Phi\egaldef \sum_{k=0}^{\infty}\phi_{k}<\infty\enspace .
\end{equation}
\begin{equation}\label{cond.ND}
\hyptag{ND}\exists \lambda<1 \telque \forall (\omega,a)\in A^{-\N}\times A,\qquad P(a|\omega)\leq \lambda\enspace .
\end{equation}
Let $d_{n}\leq n-1$, $q_{n}\leq (n-1)/4$, $t\in\set{n-1,n}$ and assume that $r_{n}\egaldef d_{n}+q_{n}\leq t$. Let $\Lst=4\paren{\Phi+\frac1{\lambda-1}}$. There exists an event $\Omega_{coup}$ satisfying $\mu\set{\Omega_{coup}^{c}}\leq 2n^{2}\beta_{q_{n}}$ such that, for all $y>0$ and all $\omega\in A^{(d_{n}+1)}$,
\begin{equation}\label{eq.conc.bas}
\P\set{\absj{\muh_{t}(\omega)-\mu(\omega)}>\Lst\sqrt{\frac{\mu(\omega)y}{n-\absj{\omega}+1}}+\frac {(d_{n}+r_{n})y}{n-\absj{\omega}+1}\cap\Omega_{coup}}\leq 4e^{-y}\enspace .
\end{equation}
\end{theorem}

\begin{remark}
An important example of application is the case of geometrically mixing processes, i.e., when the following assumption holds.
\begin{equation}\label{cond.geo.bet.mix}
\hyptag{GMC}\exists (L_{mix},\gamma_{mix})\in (\R_{+}^{*})^{2},\telque \forall k\in \N,\;\beta_{k}\leq L_{mix}e^{-\gamma_{mix}k}\enspace .
\end{equation}
Then we can choose $d_{n}=\log n$, $q_{n}=O(\log n)$ in Theorem \ref{theo.conc.bas} and we obtain that geometrically mixing processes satisfy \eqref{hyp.conc} with $\rho_{n}=(n-\log n)^{-1}\Lst^{2}$, $\varrho_{n}=O(n^{-1}\log n)$, $\varphi_{n}=n^{-2}$ and $\Omega_{conc}=\Omega_{coup}$.
\end{remark}
A immediate consequence of the previous remark is that the following corollary of Theorem \ref{theo.oracle.gal} holds.
\begin{corollary}\label{coro.BIC.oracle}
Let $(X_n)_{n\in\Z}$ be a $\phi$-mixing process satisfying \eqref{cond.phi.mix}, \eqref{cond.ND} and \eqref{cond.geo.bet.mix}. Let $d_{n}=\log(n)$, let $\pi$ be the uniform probability measure on $\Tcal_{d_n}$ and let
\begin{align*}
\Fcalst^{(n)}&=\set{\omega\in\Fcal,\telque \forall a\in A,\;(n-\absj{\omega}+1)\muh_{n}(\omega a)\in\set{0}\cup\left[(\ln n)^{4},+\infty\right[}\enspace.
\end{align*}
Let $L>18+81\ph_{\min}^{-1}$ and let $\tauh$ be the penalized estimator defined in \eqref{def.est.pen}, with
\begin{equation}\label{cond.pen.BIC}
\forall \tau\in \Fcalst^{(n)},\;\pen(\tau)\geq L\Lst^{2}|A|N(\tau)\frac{\ln n}n\enspace .
\end{equation}
There exist $n_{o}$ and a constant $\Cst$ such that, for all $n\geq n_{o}$, we have
\[
\P\set{\forall \tau\in \Fcalst^{(n)},\qquad \Cst K_{\mu}(P,\Ph_{\tauh})\leq K_{\mu}(P,\olP_{\tau})+\pen(\tau)}\geq 1-\frac2{n^{2}}\enspace .
\]
\end{corollary}

\begin{remark}
This corollary shows that BIC estimators have oracle properties, provided that the constant $c$ is sufficiently large. The drawback of this result is that the constant $\Lst^{2}$ is unknown in practice. We recommend to use the slope algorithm to overcome this problem. We will present some simulations to emphasize the advantages of this approach.
\end{remark}

\section{Simulation Study}\label{sect.simu}
In this section, we illustrate our theoretical results by simulation experiments in the family of renewal processes.
A renewal process is defined here as a binary valued process ($A=\set{0,1}$) for which the distances of successives occurrences of symbol $1$ are independent, identically distributed variables. In our simulations, the renewal distribution was Poisson with parameter $3$. The models we considered were all renewal processes with renewal times bounded by $K_o=14$. Their context trees are the subtrees of  $\tau=\set{10^{k},\;k=0,\ldots K_{o}}\cup\set{0^{K_{o}+1}}$. 
In this experiment, we used a sample size of $n=500$.
%For every $\omega\in \tau$, we chose the values $P(1|\omega)$ independent at random, uniformly distributed on $[0,1]$. We observe sequences of size $n$ and we repeat $N$ times the experiments.

\subsection{Bias and variance of the risk}
\begin{figure}%[!ht]
\begin{center}
 \includegraphics[width=12cm]{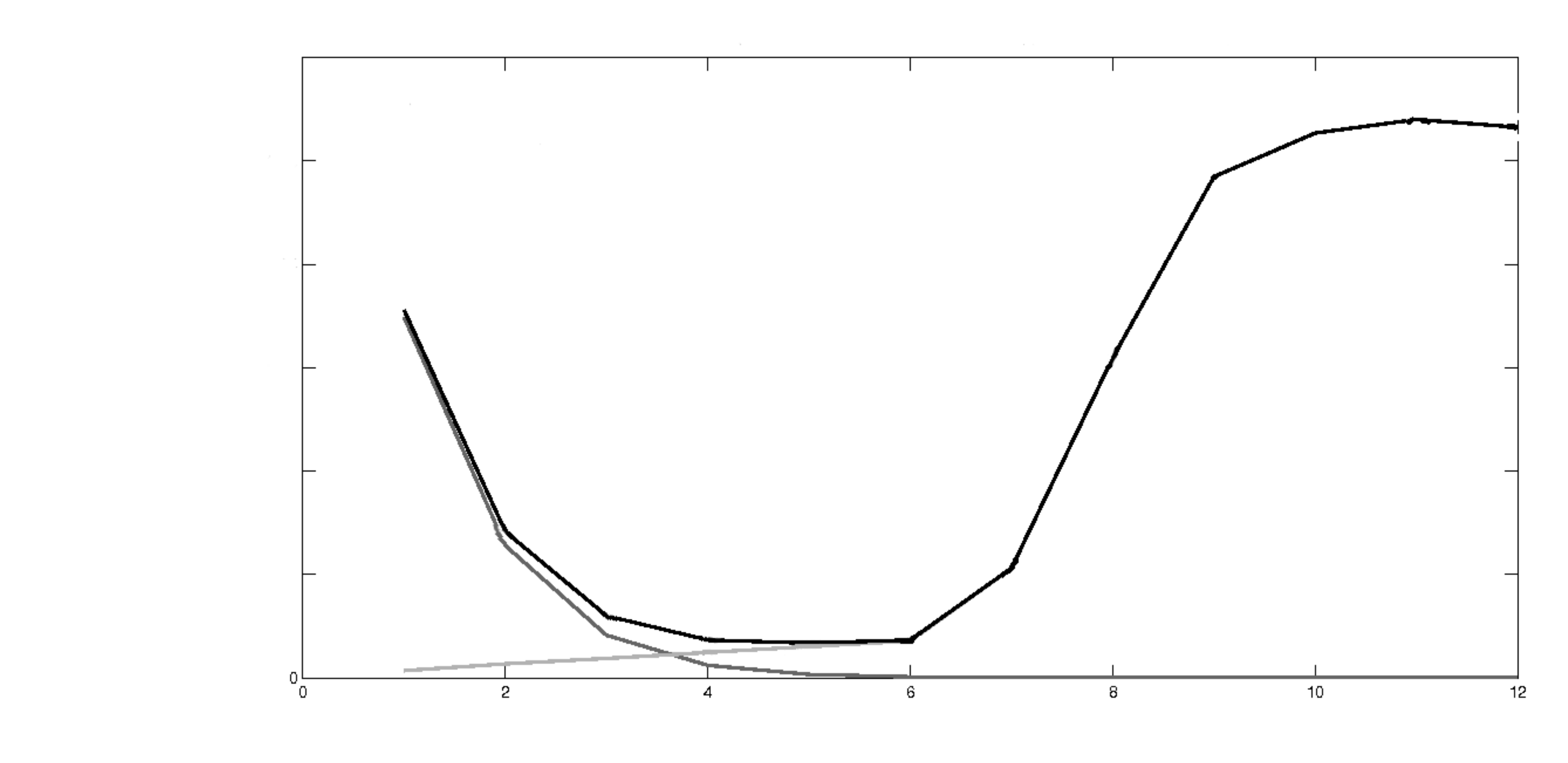}
 \caption{\label{fig.bias.variance} \small Bias (dark gray) and variance (light gray) part of the risk (black). Observe that the oracle has a size of order $5$, much smaller than the exact tree.} 
 \end{center}
\end{figure}

Figure \ref{fig.bias.variance} shows the bias and variance terms of the risk of the trees $\tau_{k_{o}}=\set{10^{k},\;k=0,\ldots k_{o}}\cup\set{0^{k_{o}+1}}$ in the previous model as a function of $k_o$. The bias can be computed easily; the variance part is estimated by a Monte-Carlo method over $N=10000$ experiments. 

\subsection{The slope phenomenon}
\begin{figure}%[!ht]
\begin{center}
 \includegraphics[width=12cm]{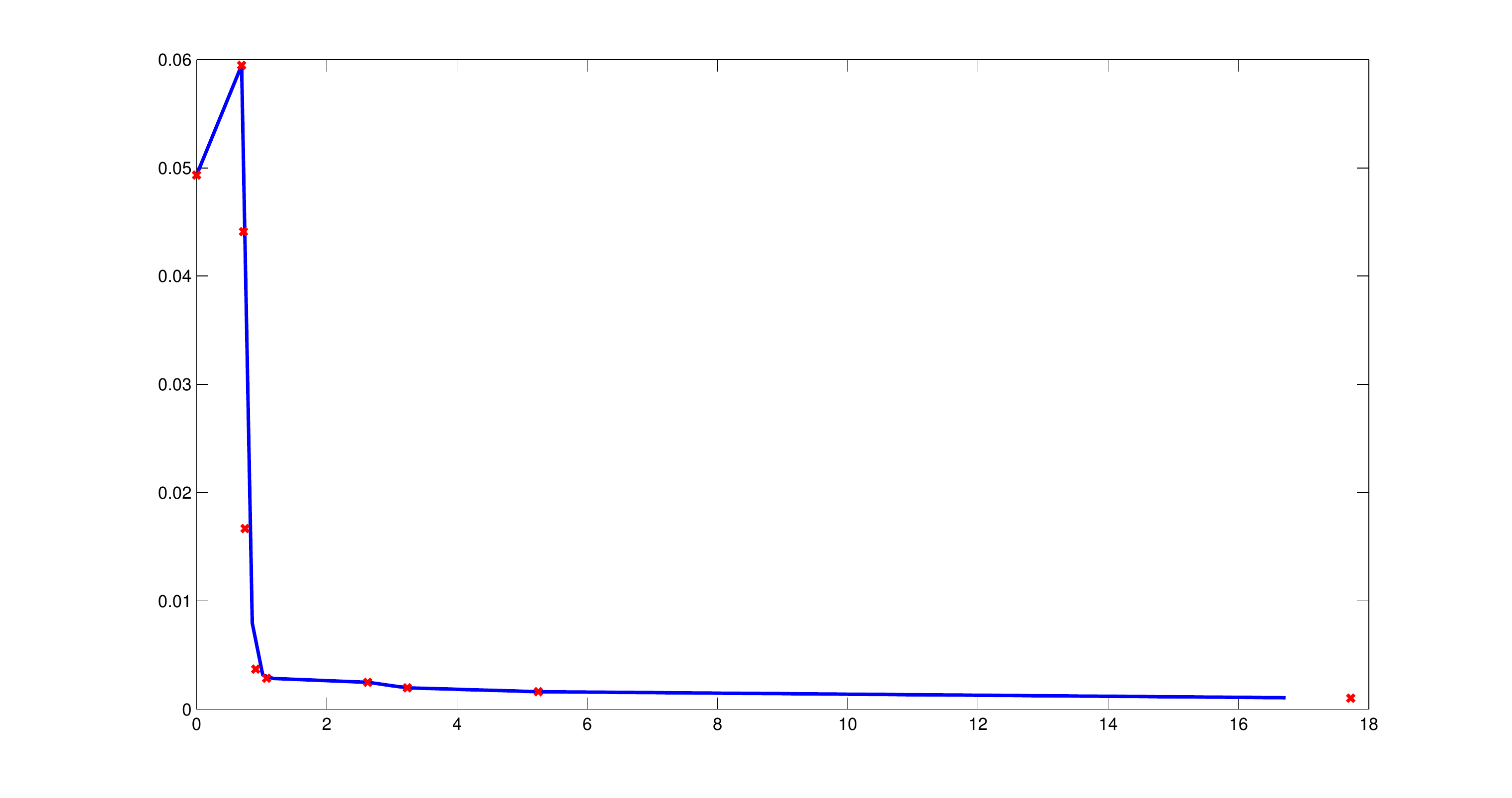}
 \caption{\label{fig.slope.heuristic} \small Example of slope heuristic. Observe the sudden change in behavior around the minimal penalty.} 
 \end{center}
\end{figure}

We illustrate the slope phenomenon. The measure of complexity is $K_{\muh}(\Ph_{\tau},\Ph^{BS}_{\tau})$, where  $\Ph^{BS}_{\tau}$ is a bootstrap estimator of $\Ph_{\tau}$ and we plot the complexity of the tree selected by minimization of the criterion
\[
\crit(\tau)=\sum_{\omega\in \tau}\muh_{n-1}(\omega)\sum_{a\in A}\Ph_{\tau}(a|\omega)\ln\paren{\frac1{\Ph_{\tau}(a|\omega)}}+cK_{\muh}(\Ph_{\tau},\Ph^{BS}_{\tau})\enspace ,
\]
for the positive constants $c$. We clearly see that when $c$ is smaller than $1$ the complexity is the largest possible and this is the content of Theorem \ref{theo.min.pen}. We also observe that when $c$ is slightly larger than $1$ there is a sudden decrease in the complexity, which is the content of Theorem \ref{theo.opt.pen}. The result are shown in Figure \ref{fig.slope.heuristic}.
For these small values of $c$, very large models are chosen, and the bootstrap estimation of their complexity is not reliable: this explains the absence of monotonicity in the left-most part of the graph, as well as in the right-most part of Figure~\ref{fig.variances}.

\subsection{Slope algorithm}
In this section, we show the performances of the slope algorithm. We take $n=500$ and $N=10000$. We use the following penalization procedures.\\

{\bf Method 1 : BIC.} The penalty is equal to the BIC penalty with $c=1/2$.\\

{\bf Method 2 : BIC+Slope.} The penalty term is equal to the BIC penalty and the constant $c$ is computed with the slope algorithm {\bf SA1-SA2-SA3} of section \ref{sect.slope.algo}. In the step {\bf SA2}, we choose for $L_{\min}$ the constant minimizing a discrete derivative of the function $L\flapp \pen_{L}(\tauh(L))$, where $\tauh(L)$ is the tree selected by $L|A|N(\tau)(\ln n)/n$.\\

{\bf Method 3 : Resampling.} For all words $\omega$, the conditional probabilities are estimated by a bootstrap method and, following Efron's heuristic (see \cite{Ef79}), $K_{\mu_{\tau}}(P_{\tau},\Ph_{\tau})$ is then estimated by the quantity $K_{\muh}(\Ph_{\tau},\Ph_{\tau}^{BS})$ and, following Theorem \ref{theo.opt.pen} the penalty is taken equal to $2K_{\muh}(\Ph_{\tau},\Ph_{\tau}^{BS})$.\\

{\bf Method 4 : Resampling+Slope.} The penalty term is  $LK_{\muh}(\Ph_{\tau},\Ph_{\tau}^{BS})$, where the constant $L$ is evaluated by the slope algorithm, with the complexity $K_{\muh}(\Ph_{\tau},\Ph_{\tau}^{BS})$. Step {\bf SA2} of the slope algorithm is evaluated in the same way as in Method 2.\\

\begin{figure}%[!ht]
\begin{center}
 \includegraphics[width=12cm]{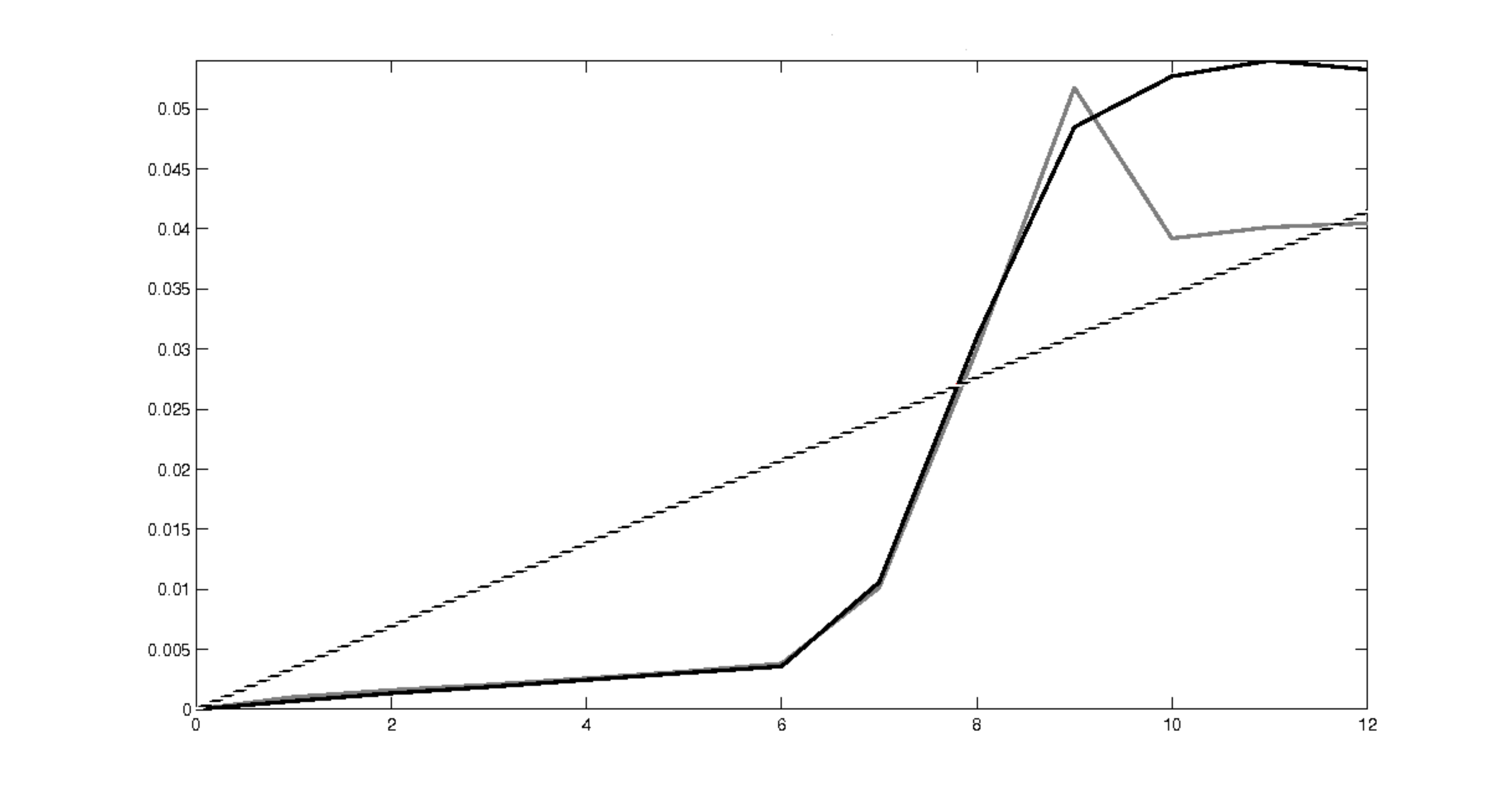}
 \caption{\label{fig.variances} \small Variance term (black), bootstrap estimator (gray) and BIC penalty with c=1/2 (dotted).} 
 \end{center}
\end{figure}

The motivation to use resampling methods comes from the fact that the variance term is better estimated than with the BIC penalty, as shown by Figure \ref{fig.variances}.

\begin{figure}%[!ht]
\begin{center}
 \includegraphics[width=12cm]{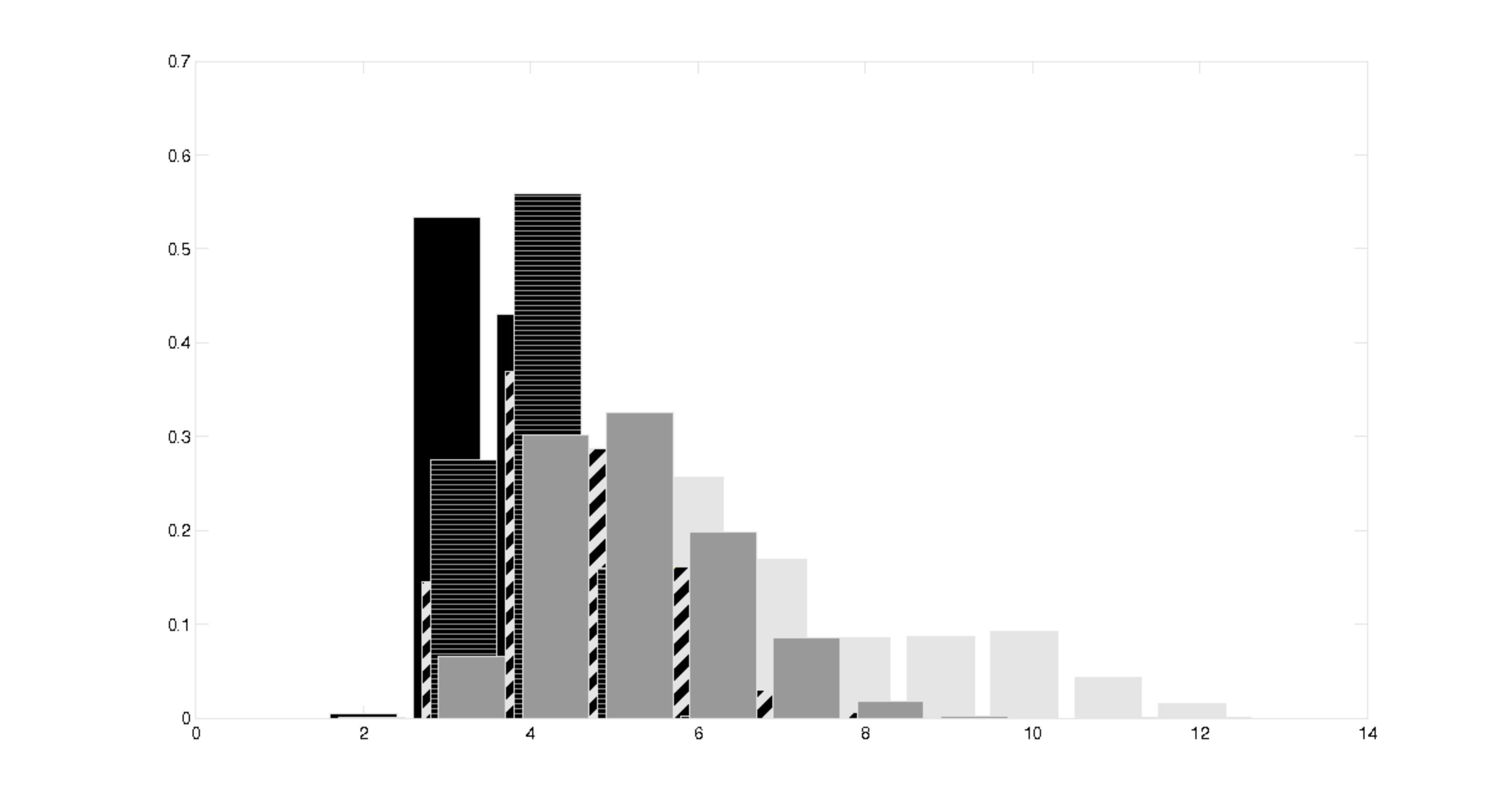}
 \caption{\label{fig.comp.methods} \small Histograms of the models selected by methods 1--4, respectively filled in black, slanted hatch, horizontal hatch and gray; models selected by AIC are depicted in light gray. Observe that the models selected by BIC penalties are smaller than those obtained with resampling penalties, both with or without the slope algorithm. As expected, for such a small sample the selected estimator is always significantly smaller than the actual context tree of the source.} 
 \end{center}
\end{figure}

Figure \ref{fig.comp.methods} presents histograms of the models selected by methods 1--4, for $n=1000$, $N=10000$.

Finally, in order to illustrate the oracle properties of the selected estimators, we compute in the table \ref{table.compare} the values of the ratio 
\[
\frac{K_{\mu}(P,\Pt_{\tauh})}{\inf_{\tau\in\Fcal_{n-1}}K(P,\Pt_{\tau})}
\]
for the different methods. We give the mean value over $N=10000$ experiments of the risk ratio and the standard deviation is also indicated.
\begin{table}[!ht]
\begin{center}
    \caption{Comparative table of oracle ratios}
    \begin{tabular}{ | p{1.3cm} || p{2.2cm} | p{2.2cm} | p{2.2cm} | p{2.2cm} |} 
    %\hline
   Method & BIC & BIC+Slope & Resampling & Res+Slope  \\ \hline
   risk ratio & 1.5245 (0.8568) & 1.2665    (0.5657) & 1.2751    (0.3230) & 1.6707    (1.9702) \\ \hline    
    \end{tabular}
    \label{table.compare}
\end{center}
\end{table}

It is interesting to remark that the oracle performances of the BIC estimator are improved by the slope algorithm. On the other hand, resampling estimators do not seem to improve significantly the results. Moreover, the slope algorithm combined with this penalization method give the worst results here. As the computational cost of resampling methods is quite heavy, we do not recommend to use it in practice. On the other hand, the slope algorithm does not add a significant computational cost and can be used to choose the leading constant.\\

Note that, in general, Methods 3 and 4 involve a minimization problem \eqref{def.est.pen} that is not computationally tractable.
Here, the particular structure of renewal processes allowed us to consider all the models; in more general settings, we propose to proceed in two steps: first, select a set of trees defined by the image of $c\flapp \tauh_{BIC}(c)$ $c>0$, where, for all $c>0$, $\tauh_{BIC}(c)$ is the tree selected by the BIC penalty with constant $c$. Then, select among those trees with the proposed resampling methods.

\section{Conclusion}

We developed an oracle approach for context tree selection with K\"ullback loss.  %The resulting estimator has a natural interpretation in information theory as it yields a code with minimal redundancy. 
Our presentation emphasizes the central role of concentration inequalities for sequences of words. We proved that such concentration inequalities hold for geometrically $\phi$-mixing sequences. We obtained as a corollary of our general approach some oracle properties of the BIC-like estimators in this framework.

We also provided both numerical and theoretical justification for the use of the  slope heuristic in this problem in order to calibrate the leading constant in the penalties. This provides in particular an answer for the practical choice of the leading constant in the BIC-like penalty. Actually, \cite{CT06} proved the consistency of the BIC estimators for any value of $c$. On the other hand, \cite{GGGL} proved that, for any finite value of $n$, the set of trees selected by the BIC-like penalties $(\pen_{c}(\tau))_{c>0}$ is the set of champions, where, for any $k\leq \log n$, the champion of size $k$ is the one maximizing the log-likelihood among those trees $\tau$ such that $N(\tau)\leq k$.

There is a growing interest for the slope heuristic, see for example \cite{BM07, AM08, Le09, Le10, Sa11, LT11}. However, the theoretical analysis of this method is still in its beginning and our results are a significative contribution. In particular, we provide, up to our knowledge, the first proof of the relevance of the slope heuristic in a discrete non i.i.d framework.

Our results also emphazise the interest of the oracle approach, compared to the identification approach. In fact, a large part of the interest of context tree models lies in the fact that \emph{any} stationary ergodic source can be approached, in the K\"ullback information distance, by context tree sources: hence, the use of these models is not restricted to cases when the true source belongs to one of them. Besides, even if it is finite, the true source's context tree is likely not to be the best model to use for small samples, as illustrated in our simulation study. In fact, we showed that the BIC estimator presents nice oracle properties, and that it can be further improved by choosing the leading constant in the penalty adaptively. 
This result justifies the use of context tree models in practical applications much more than the consistency properties that are usually mentioned.

An important question related to the oracle approach is to obtain upper bounds for the risk of the selected estimator. We showed in Section \ref{sec.cont.rates} that such bounds can be obtained with continuity rates. Actually, these continuity rates provide upper bounds for the bias term and yield good mixing properties, so that we can use the upper bounds of the variance term obtained in the mixing case.

The $\phi$-mixing properties assumed in Section \ref{sect.mix} are somewhat restrictive. It would be interesting to work with weaker assumptions, for example, with weaker mixing coefficients, as $\tilde{\phi}$ (see \cite{DP05} for a definition). These mixing coefficients are sufficient to generalize some results of model selection (see \cite{Le08,Le10} for example). Another interesting problem would be to look for natural mixing properties of context tree sources. As mentioned, we proved such properties in Section \ref{sec.cont.rates} using a theorem of \cite{CFF02}. This last theorem was obtained as a consequence of the existence of a constructive perfect simulation scheme for the chains. New perfect simulation schemes have been developed recently \cite{garivier11PWCT}, using less restrictive assumptions on the chains. It would be interesting to see what mixing-properties can be deduced from these new constructions.

\subsection*{Acknowledgements} 

We would like to thank gratefully Roberto I. Oliveira who pointed out the results of Section \ref{section.continuity.mixing}. We also want to thank Antonio Galves for many discussions and fruitful advices during the redaction of the paper.

\bibliographystyle{imsart-number}
\bibliography{bibliolerasle}

\appendix
\section{Proofs in the general case}
\subsection{Proof of Proposition \ref{prop.typ.gal}}
Let $(\omega,a)\in \Fcalst(\delta)\times A$. By definition, if $\mu(\omega a)\neq 0$,
\[
\sqrt{\rho_{n}\mu(\omega a)\ln\paren{\frac1{\pi(\omega a)\delta}}}+\varrho_{n}\ln\paren{\frac1{\pi(\omega a)\delta}}\leq \paren{\sqrt{\frac1{\Lambda_{n}^{(1)}}}+\frac1{\Lambda_{n}^{(2)}}}\mu(\omega a)\enspace .
\]
Moreover,
\begin{align*}
\mu(\omega)\geq \mu(\omega a)&\geq \ln\paren{\frac1{\pi(\omega a)\delta}}\paren{\Lambda_{n}^{(1)}\rho_{n}\vee \Lambda_{n}^{(2)}\varrho_{n}}\\
&\geq \ln\paren{\frac1{\pi(\omega)\delta}}\paren{\Lambda_{n}^{(1)}\rho_{n}\vee \Lambda_{n}^{(2)}\varrho_{n}}\enspace .
\end{align*}
Hence,
\[
\sqrt{\rho_{n}\mu(\omega)\ln\paren{\frac1{\pi(\omega)\delta}}}+\varrho_{n}\ln\paren{\frac1{\pi(\omega)\delta}}\leq \paren{\sqrt{\frac1{\Lambda_{n}^{(1)}}}+\frac1{\Lambda_{n}^{(2)}}}\mu(\omega)\enspace .
\]
Thus, on $\Omega_{good}$, all the words in $\Fcalst(\delta)$ belong to $\Tbf_{typ}(\eta)$, with
\[
\eta=\paren{\sqrt{\frac1{\Lambda_{n}^{(1)}}}+\frac1{\Lambda_{n}^{(2)}}}\enspace .
\]
Let now $(\omega,a)\in \Fcalst^{(n)}(\delta)\times A$. On $\Omega_{good}$, if $\muh_{n}(\omega a)\neq 0$, we have
\[
\muh_{n}(\omega a)\leq \paren{\sqrt{\mu(\omega a)}+\sqrt{\rho_{n}\ln\paren{\frac1{\pi(\omega)\delta}}}}^{2}+(\varrho_{n}-\rho_{n})\ln\paren{\frac1{\pi(\omega)\delta}}\enspace .
\]
Hence, by definition of $\Fcalst^{(n)}(\delta)$,
\begin{equation}\label{eq.F.sub.Fh}
\paren{\sqrt{1-\frac1{\Lambda_{n}^{(2)}}}-\sqrt{\frac1{\Lambda_{n}^{(1)}}}}^{2}\muh_{n}(\omega a)\leq \mu(\omega a)\enspace .
\end{equation}
As a consequence, $\omega\in\Fcalst(\delta)$ for all $n$ such that 
\[
\paren{\sqrt{1-\frac1{\Lambda_{n}^{(2)}}}-\sqrt{\frac1{\Lambda_{n}^{(1)}}}}^{2}\geq \frac12\enspace .
\]
Let finally $(\omega,a)\in \Fcalst^{(2)}(\delta)\times A$. On $\Omega_{good}$, if $\mu(\omega a)\neq 0$, we have
\[
\mu(\omega a)\paren{1-\sqrt{\frac1{\Lambda_{n}^{(1)}}}-\frac1{\Lambda_{n}^{(2)}}}\leq \muh_{n}(\omega a)\enspace .
\]
Hence, by definition of $\Fcalst^{(2)}(\delta)$, $\omega\in\Fcalst^{(n)}(\delta)$ for all $n$ such that 
\[
\sqrt{\frac1{\Lambda_{n}^{(1)}}}-\frac1{\Lambda_{n}^{(2)}}\leq \frac12\enspace .
\]

\subsection{Proof of Theorem \ref{theo.oracle.gal}}
Let $\tau\subset \Fcalst^{(n)}(\delta)$. From Lemma \ref{lem.dec.oracle} and the definition of $\tauh$, we have
\begin{align}
\notag K_{\mu}(P,\Ph_{\tauh})\leq &K_{\mu}(P,\olP_{\tau})-K_{\muh}(\Ph_{\tau},P_{\tau})+\pen(\tau)\\
\label{eq.dec.risk}&+\paren{K_{\mu_{\tauh}}(P_{\tauh},\Ph_{\tauh})+K_{\muh}(\Ph_{\tauh},P_{\tauh})-\pen(\tauh)}+L(P_{\tau})-L(P_{\tauh})\enspace .
\end{align}
It comes from Proposition \ref{prop.typ.gal} that, on $\Omega_{good}$, all the words in $\tau\cup\tauh$ belong to $\Tbf_{typ}(\eta)$ for some $\eta=O\paren{\sqrt{\frac1{\Lambda_{n}^{(1)}}}\vee \frac1{\Lambda_{n}^{(2)}}}$. Therefore, Lemma \ref{lem.termes.croises.2} gives, for any $\epsilon>0$,
\begin{align*}
L(P_{\tau})-L(P_{\tauh})
\leq & (\epsilon+O(\eta)) \paren{K_{\mu}(P,\olP_{\tau})+K_{\mu}(P,\olP_{\tauh})}\\
&+(1+O(\eta))\frac{\Lst^{(\tau,\tauh)}} {\epsilon}\paren{K_{\mu_{\tau}}(P_{\tau},\Ph_{\tau})+K_{\mu_{\tauh}} (P_{\tauh},\Ph_{\tauh})}\enspace .
\end{align*}
Assume now that $\eta\leq 1/3$ and let $\epsilon=1/2$. By typicality, it holds that $\Lst^{(\tau,\tauh)}\leq 3\ph_{\min}^{-1}$, and hence, from \eqref{eq.dec.risk}, we deduce that $(1/2-O(\eta)) K_{\mu}(P,\Ph_{\tauh})$ is upper-bounded by
\begin{align}
\label{eq.dec.risk.2}&\left( \frac{3}{2}+O(\eta)\right)K_{\mu}(P,\olP_{\tau})+(1+O(\eta))\frac{6}{\ph_{\min}}K_{\mu_{\tau}}(P_{\tau},\Ph_{\tau})-K_{\muh}(\Ph_{\tau},P_{\tau})\\
&+\pen(\tau)+(1+O(\eta))\paren{1+\frac{6}{\ph_{\min}}}K_{\mu_{\tauh}}(P_{\tauh},\Ph_{\tauh})+K_{\muh}(\Ph_{\tauh},P_{\tauh})-\pen(\tauh)\; .\nonumber
\end{align}
In addition, from Lemma \ref{lem.cons.typ}, there exists $\eta^{\prime}=O\paren{(\Lambda_{n}^{(1)}\wedge\Lambda_{n}^{(2)})^{-1/2}}$ such that, for all $6+18\ph_{\min}^{-1}<L^{\prime}<L$, on $\Omega_{good}$
\begin{align*}
\paren{1+O(\eta)}&\paren{1+\frac{6}{\ph_{\min}}}K_{\mu_{\tauh}}(P_{\tauh},\Ph_{\tauh})+K_{\muh}(\Ph_{\tauh},P_{\tauh})\\
\leq& (L^{\prime}+\eta^{\prime})\paren{\sqrt{\rho_{n}}+\sqrt{ \frac{\varrho_{n}}{\Lambda_{n}^{(2)}}}}^{2}\sum_{(\omega,a)\in \tauh\times A}\ln\paren{\frac1{\pi(\omega a)\delta}}\;.
%\\&+\paren{\paren{(1+O(\eta))\paren{2+\frac6{\ph_{\min}}}-\frac{L^{\prime}}{3}}\wedge 0}K_{\mu_{\tauh}}(P_{\tauh},\Ph_{\tauh})\enspace .
\end{align*}
For $n$ sufficiently large, we have $L^{\prime}+\eta\leq L$, hence
\begin{align}\label{eq.main.term}
\Big(1+&O(\eta)\Big)\paren{1+\frac{6}{\ph_{\min}}}K_{\mu_{\tauh}}(P_{\tauh},\Ph_{\tauh})+K_{\muh}(\Ph_{\tauh},P_{\tauh})%\nonumber\\&
\leq \pen(\tauh)\;.%+\paren{\paren{(1+O(\eta))\paren{2+\frac6{\ph_{\min}}}-\frac{L^{\prime}}{3}}\wedge 0}K_{\mu_{\tauh}}(P_{\tauh},\Ph_{\tauh})\enspace .
\end{align}
We conclude the proof, plugging \eqref{eq.main.term} in \eqref{eq.dec.risk.2}.

\subsection{Proof of Theorem \ref{theo.min.pen}}
Thanks to Proposition \ref{prop.typ.gal}, there exists $n_{o}$ such that for $n\geq n_{o}$, on $\Omega_{good}$, $\Fcalst^{(n)}(\delta)\subset \Fcalst(\delta)$. Let $\tau_{-}$ be any element in $\Fcalst(\delta)$. $\tauh$ minimizes over $\Fcalst^{(n)}(\delta)$ the following criterion:
\begin{multline*}
\crit_{\tau_{-}}(\tau)\egaldef \sum_{(\omega,a)\in\tau\times A}\muh_{n}(\omega a)\ln\paren{\frac1{\Ph_{\tau}(a|\omega)}}\\+\int_{A^{-\N}\times A}d\mu(\omega a)\ln(P(a|\omega))-L(P_{\tau_{-}})+\pen(\tau)\enspace .
\end{multline*}
Thanks to Lemma \ref{lem.dec.oracle}, we have
\[
\crit_{\tau_{-}}(\tau)=K_{\mu}(P,\olP_{\tau})-K_{\muh}(\Ph_{\tau},P_{\tau})+L(P_{\tau})-L(P_{\tau_{-}})+\pen(\tau)\enspace .
\]
Thanks to Proposition \ref{prop.typ.gal}, there exists $\eta_{1}=O\paren{\sqrt{\frac1{\Lambda_{n}^{(1)}}}\vee\frac1{\Lambda_{n}^{(2)}}}$ such that $\Fcalst(\delta)\subset \Tbf_{typ}(\eta_{1})$. Hence, from \eqref{eq.cons.typ.3} in Lemma \ref{lem.cons.typ.gal}, for all $\tau\in \Fcalst(\delta)$, on $\Omega_{good}$
\[
\absj{K_{\muh}(\Ph_{\tau},P_{\tau})-K_{\mu_{\tau}}(P_{\tau},\Ph_{\tau})}\leq O(\eta_{1})K_{\mu_{\tau}}(P_{\tau},\Ph_{\tau})\enspace .
\]
In addition, from Lemma \ref{lem.termes.croises.2}, for all $\tau\neq \tau_{-}$, on $\Omega_{good}$, we have 
\[|L(P_{\tau})-L(P_{\tau_{-}})|\leq A(\tau_{-}, \tau)\;,\]
where
\begin{align}
 \notag A(\tau_{-}, \tau)&\egaldef
 (2+\eta_{1})\sqrt{\Lst^{(\tau_{-},\tau)}} \times \\
\label{def.Atau}&\sqrt{\paren{K_{\mu}(P,\olP_{\tau_{-}})+K_{\mu}(P,\olP_{\tau})}\paren{K_{\mu_{\tau}}(P_{\tau},\Ph_{\tau})+K_{\mu_{\tau_{-}}}(P_{\tau_{-}},\Ph_{\tau_{-}})}}\; .
\end{align}
The inequalities $\crit_{\tau_{o}}(\tauh)\leq \crit_{\tau_{o}}(\tau_{o})$ and $\crit_{\taust}(\tauh)\leq \crit_{\taust}(\taust)$ can therefore be rewritten
\begin{multline} K_{\mu}(P,\olP_{\tau_{o}})-(r-\eta_{1})K_{\mu_{\tau_{o}}}(P_{\tau_{o}},\Ph_{\tau_{o}})\\\geq K_{\mu}(P,\olP_{\tauh})-(1+\eta_{1})K_{\mu_{\tauh}}(P_{\tauh},\Ph_{\tauh})-A(\tau_o,\tauh)\label{eq.min.pen.good}\;,
\end{multline}
\begin{multline} 
K_{\mu}(P,\olP_{\taust})-(r-\eta_{1})K_{\mu_{\taust}}(P_{\taust},\Ph_{\taust})\\\geq K_{\mu}(P,\olP_{\tauh})-(1+\eta_{1})K_{\mu_{\tauh}}(P_{\tauh},\Ph_{\tauh})-A(\taust, \tauh)\label{eq.min.pen.good.2}\;.
\end{multline}
Recall that, on the event $\Omega_{hyp}$, 
\[
K_{\mu}(P,\olP_{\taust})=o\paren{ K_{\mu_{\taust}}(P_{\taust},\Ph_{\taust})}\enspace .
\]
Inequality \eqref{eq.min.pen.good.2} can then be satisfied only if one of the following condition holds.
\begin{align}
\label{cond.min.pen.1}\hyptag{C1}\exists L\egaldef L(r,\ph_{\min})&\telque K_{\mu_{\tauh}}(P_{\tauh},\Ph_{\tauh})\geq L K_{\mu_{\taust}}(P_{\taust},\Ph_{\taust})\enspace .\\
\label{cond.min.pen.2}\hyptag{C2}\exists L\egaldef L(r,\ph_{\min})&\telque K_{\mu_{\tauh}}(P_{\tauh},\Ph_{\tauh})=o\paren{ K_{\mu_{\taust}}(P_{\taust},\Ph_{\taust})}\;\\
&\nonumber\quad\quad\mbox{and}\;K_{\mu}(P,\olP_{\tauh})\geq L K_{\mu_{\taust}}(P_{\taust},\Ph_{\taust})\enspace .
\end{align}
In fact, under \eqref{cond.min.pen.2}, 
\[
K_{\mu_{\tau}}(P_{\tau},\Ph_{\tau})+K_{\mu_{\tau_{o}}}(P_{\tau_{o}},\Ph_{\tau_{o}}) = o \left(K_{\mu_{\taust}}(P_{\taust},\Ph_{\taust}) \right)=o \left( K_{\mu}(P,\olP_{\tauh}) \right)\;,
\]
hence $A(\tau_{o}, \tauh) = o \left( K_{\mu}(P,\olP_{\tauh}) \right)$. Thus inequality \eqref{eq.min.pen.good} and $K_{\mu}(P,\olP_{\tau_{o}})=o\paren{K_{\mu_{\taust}}(P_{\taust},\Ph_{\taust})}$ yield 
\[
K_{\mu}(P,\olP_{\tauh})=o \left( K_{\mu}(P,\olP_{\tauh}) \right)\enspace .
\]
This is a contradiction. Hence, condition \eqref{cond.min.pen.1} is fulfilled. By typicality, we have $p_{\min}/2\leq\ph_{\min}\leq 3p_{\min}/2$ for $n$ sufficiently large, thus $L(r,\ph_{\min})\geq L^{\prime}(r,p_{\min})$ which concludes the proof of the Theorem.
\subsection{Proof of Theorem \ref{theo.opt.pen}}
$\tauh$ minimizes over $\Fcalst^{(n)}(\delta)$ the following criterion
\begin{multline*}
\crit(\tau)\egaldef \sum_{(\omega,a)\in\tau\times A}\muh_{n}(\omega a)\ln\paren{\frac1{\Ph_{\tau}(a|\omega)}}\\
+\int_{A^{-\N}\times A}d\mu(\omega a)\ln(P(a|\omega))-L(P_{\tau_{o}})+\pen(\tau)\enspace .
\end{multline*}
Thanks to Lemma \ref{lem.dec.oracle}, we have
\[
\crit(\tau)=K_{\mu}(P,\olP_{\tau})-K_{\muh}(\Ph_{\tau},P_{\tau})+L(P_{\tau})-L(P_{\tau_{o}})+\pen(\tau)\enspace .
\]
Thanks to Proposition \ref{prop.typ.gal}, there exists $\eta_{1}=O\paren{\sqrt{\frac1{\Lambda_{n}^{(1)}}}\vee\frac1{\Lambda_{n}^{(2)}}}$ such that $\Fcalst(\delta)\subset \Tbf_{typ}(\eta_{1})$. Hence, from \eqref{eq.cons.typ.3} in Lemma \ref{lem.cons.typ.gal}, for all $\tau\in \Fcalst(\delta)$, on $\Omega_{good}$
\[
\absj{K_{\muh}(\Ph_{\tau},P_{\tau})-K_{\mu_{\tau}}(P_{\tau},\Ph_{\tau})}\leq O(\eta_{1})K_{\mu_{\tau}}(P_{\tau},\Ph_{\tau})\enspace .
\]
In addition, from Lemma \ref{lem.termes.croises.2}, for all $\tau\neq \tau_{o}$, on $\Omega_{good}$, we have 
\[|L(P_{\tau})-L(P_{\tau_{o}})|\leq A(\tau_o,\tau)\;,\]
where $A(\tau,\tau')$ is defined in \eqref{def.Atau}.%
The inequalities $\crit(\tauh)\leq \crit(\tau_{o})$ can therefore be rewritten
\begin{multline}
\label{eq.opt.pen.1} K_{\mu}(P,\olP_{\tau_{o}})+(r_{2}+O(\eta_{1}))K_{\mu_{\tau_{o}}}(P_{\tau_{o}},\Ph_{\tau_{o}})\\\geq K_{\mu}(P,\olP_{\tauh})+(r_{1}-O(\eta_{1}))K_{\mu_{\tauh}}(P_{\tauh},\Ph_{\tauh})
-A(\tau_o, \tauh)\;.
\end{multline}
In $A(\tau_o, \tauh)$, all the terms are, on $\Omega_{hyp}$, $o(K_{\mu_{\taust}}(P_{\taust},\Ph_{\taust}))$, therefore, \eqref{eq.opt.pen.res.1} follows from \eqref{eq.opt.pen.1}. Moreover, using repeatedly the inequalities, valid for any $a>0$, $b>0$, $\epsilon>0$, 
\[
\sqrt{a+b}\leq \sqrt{a}+\sqrt{b}\;\mbox{and}\;2\sqrt{ab}\leq \epsilon a+\frac b{\epsilon}\enspace,
\]
we obtain, for any $\epsilon_{1}>0$, $\epsilon_{2}>0$,
\begin{align*}
2\sqrt{\Lst^{(\tau_{o},\tauh)}}&\sqrt{\paren{K_{\mu}(P,\olP_{\tau_{o}})+K_{\mu}(P,\olP_{\tauh})}\paren{K_{\mu_{\tauh}}(P_{\tauh},\Ph_{\tauh})+K_{\mu_{\tau_{o}}}(P_{\tau_{o}},\Ph_{\tau_{o}})}}\\
&\leq \enspace (\epsilon_{1}+\epsilon_{2}) K_{\mu}(P,\olP_{\tauh})+\paren{\epsilon_{1}+\frac{\Lst^{(\tau_{o},\tauh)}}{\epsilon_{2}}}K_{\mu_{\tauh}}(P_{\tauh},\Ph_{\tauh})\\
&+\paren{1+\frac{\Lst^{(\tau_{o},\tauh)}}{\epsilon_{1}}}K_{\mu}(P,\olP_{\tau_{o}})+\Lst^{(\tau_{o},\tauh)}\paren{1+\epsilon_{1}^{-1}}K_{\mu_{\tau_{o}}}(P_{\tau_{o}},\Ph_{\tau_{o}}).
\end{align*}
We plug this inequality in \eqref{eq.opt.pen.1}, we obtain 
\begin{multline*}
 \paren{2+\frac{\Lst^{(\tau_{o},\tauh)}}{\epsilon_{1}}+O(\eta_{1})} K_{\mu}(P,\olP_{\tau_{o}})\\+(r_{2}+\Lst^{(\tau_{o},\tauh)}\paren{1+\epsilon_{1}^{-1}}+ O(\eta_{1}))K_{\mu_{\tau_{o}}}(P_{\tau_{o}},\Ph_{\tau_{o}})\\
\geq (1-\epsilon_{1}-\epsilon_{2})K_{\mu}(P,\olP_{\tauh})+(r_{1}-\epsilon_{1}-\frac{\Lst^{(\tau_{o},\tauh)}}{\epsilon_{2}}-O(\eta_{1}))K_{\mu_{\tauh}}(P_{\tauh},\Ph_{\tauh}).
\end{multline*}

\subsection{Proof of Theorem \ref{theo.opt.pen.afaire}}
Thanks to Proposition \ref{prop.typ.gal}, there exists $n_{o}$ such that for $n\geq n_{o}$, on $\Omega_{good}$, $\Fcalst^{(n)}(\delta)\subset \Fcalst(\delta)$. $\tauh$ minimizes over $\Fcalst^{(n)}(\delta)$ the following criterion
\begin{multline*}
\crit(\tau)\egaldef \sum_{(\omega,a)\in\tau\times A}\muh_{n}(\omega a)\ln\paren{\frac1{\Ph_{\tau}(a|\omega)}}+\\\int_{A^{-\N}\times A}d\mu(\omega a)\ln(P(a|\omega))-L(P_{\tau_{o}})+\pen(\tau)\enspace .
\end{multline*}
Thanks to Lemma \ref{lem.dec.oracle}, we have
\[
\crit(\tau)=K_{\mu}(P,\olP_{\tau})-K_{\muh}(\Ph_{\tau},P_{\tau})+L(P_{\tau})-L(P_{\tau_{o}})+\pen(\tau)\enspace .
\]
Thanks to Proposition \ref{prop.typ.gal}, there exists $\eta_{1}=O\paren{\sqrt{\frac1{\Lambda_{n}^{(1)}}}\vee\frac1{\Lambda_{n}^{(2)}}}$ such that $\Fcalst(\delta)\subset \Tbf_{typ}(\eta_{1})$. Hence, from \eqref{eq.cons.typ.3} in Lemma \ref{lem.cons.typ.gal}, for all $\tau\in \Fcalst(\delta)$, on $\Omega_{good}$
\[
\absj{K_{\muh}(\Ph_{\tau},P_{\tau})-K_{\mu_{\tau}}(P_{\tau},\Ph_{\tau})}\leq O(\eta_{1})K_{\mu_{\tau}}(P_{\tau},\Ph_{\tau})\enspace .
\]
In addition, on $\Omega_{mis}$, we have
\[
L(P_{\tau})-L(P_{\tau_{o}})\leq uK_{\mu}(P,\Pt_{\tau})+vK_{\mu}(P,\Pt_{\tau_{o}})\enspace .
\]
Hence, the equation $\crit(\tauh)\leq \crit(\tau_{o})$ implies, with the conditions on the penalty
\begin{multline*}
(1-u)K_{\mu}(P,\olP_{\tau})+(1+r_{1}-u-\eta_{1})K_{\mu_{\tauh}}(P_{\tauh},\Ph_{\tauh})\\\leq (1+v)K_{\mu}(P,\olP_{\tau})+(1+r_{2}+v+\eta_{1})K_{\mu_{\tauh}}(P_{\tauh},\Ph_{\tauh})\enspace .
\end{multline*}
\section{Proofs in the mixing case}

\subsection{Proof of Theorem \ref{theo.conc.bas}}
Let us write $t=p_{n}r_{n}+u_{n}$, with $0\leq u_{n}<r_{n}$. Let us now denote, for all $k=1,\ldots ,p_{n}$, the set $I_{p_{n}+1-k}$ defined as:
\begin{itemize}
 \item $I_{p_{n}+1-k}=\set{1\vee \croch{t-d_{n}-kr_{n}+1},\ldots,t-(k-1)r_{n}}$ if $k\leq p_{n}$;
 \item $I_{p_{n}+1-k}=\set{1,\ldots,u_{n}}$ if $k=p_{n}+1$ and $u_{n}\geq d_{n}$;
 \item $I_{p_{n}+1-k}=\emptyset$ if $k=p_{n}+1$ and $u_{n}<d_{n}$.
\end{itemize}

Let $h_{1}=\un_{u_{n}\geq d_{n}}$, $k_{n}=\PEInf{(p_{n}-1+h_{1})/2}$, $\ell_{n}=\PEInf{(p_{n}-2+h_{1})/2}$. 
We apply Lemma \ref{lem.Vi97} to the process $(X_{i})_{i\in \Z}$ and to the sets $(J_{k})_{k=1,\ldots,k_{n}}$ and $(J_{k}^{\prime})_{k=1,\ldots,\ell_{n}}$ where, for all $k=0,\ldots,k_{n}$, $J_{k}=I_{1-h_{1}+2k}$ and, for all $k=0,\ldots,\ell_{n}$, $J_{k}^{\prime}=I_{2-h_{1}+2k}$.
We obtain the random variables $(Y_{i})_{i\in \cup_{k=0}^{k_{n}}J_{k}}$ and $(Y_{i}^{\prime})_{i\in \cup_{k=0}^{\ell_{n}}J_{k}^{\prime}}$ such that,
\begin{enumerate}
\item for all $k=0,\ldots,k_{n}$, $(Y_{i})_{i\in J_{k}}$ has the same distribution as $(X_{i})_{i\in J_{k}}$ and, for all $k=0,\ldots,\ell_{n}$, $(Y_{i}^{\prime})_{i\in J_{k}^{\prime}}$ has the same distribution as $(X_{i})_{i\in J_{k}^{\prime}}$
\item for all $k=1,\ldots,k_{n}$, $(Y_{i})_{i\in J_{k}}$ is independent of $(X_{i},Y_i)_{i\in \cup_{t\leq k-1}J_{t}}$ and, for all $k=1,\ldots,\ell_{n}$, $(Y_{i}^{\prime})_{i\in J_{k}^{\prime}}$ is independent of $(X_{i},Y_{i})_{i\in \cup_{t\leq k-1}J_{t}^{\prime}}$,
\item for every $k=0,\ldots,k_{n}$, $\P\set{(X_{i})_{i\in J_{k}}\neq (Y_{i})_{i\in J_{k}}}\leq (r_{n}+d_{n})\beta_{q_{n}}$ and, for every $k=0,\ldots,\ell_{n}$, $\P\set{(X_{i})_{i\in J_{k}^{\prime}}\neq (Y_{i}^{\prime})_{i\in J_{k}^{\prime}}}\leq (r_{n}+d_{n})\beta_{q_{n}}$.
\end{enumerate}
Let $\Omega_{coup}$ be the following set
\begin{multline}\label{def.Omega.coup}
\Omega_{coup}=\Bigg\{\forall k=0,\ldots,k_{n},\; (X_{i})_{i\in J_{k}}= (Y_{i})_{i\in J_{k}}\\\;\mbox{and}\;\forall k=0,\ldots,\ell_{n},\;(X_{i})_{i\in J_{k}^{\prime}}= (Y_{i}^{\prime})_{i\in J_{k}^{\prime}}\Bigg\}\enspace .
\end{multline}
It comes from point 3 that $\P\set{\Omega_{coup}^{c}}\leq (k_{n}+\ell_{n}+2)(r_{n}+d_{n})\beta_{q_{n}}\leq 2n^{2}\beta_{q_{n}}$. Let now $\omega\in A^{(d_{n}+1)}$ such that $\mu(\omega)\neq 0$. For every $k\leq t$, let 
\[
Z_{k}=\frac{\un_{X_{k-\absj{\omega}+1}^{k}=\omega}-\mu(\omega)}{\sqrt{\mu(\omega)}}\enspace .
\]
Let also, if $\exists k\in\set{0,\ldots,\ell_{n}}\telque \set{i-\absj{\omega}+1,\ldots, i}\subset J_{k}^{\prime}$
\[
Z_{i}^{\prime}=\frac{\un_{(Y^{\prime})_{i-\absj{\omega}+1}^{i}=\omega}-\mu(\omega)}{\sqrt{\mu(\omega)}}\;,
\]
otherwise, let 
\[
Z_{i}^{\prime}=\frac{\un_{Y_{i-\absj{\omega}+1}^{i}=\omega}-\mu(\omega)}{\sqrt{\mu(\omega)}} \;.
\]
For $k=1-h_{1},\ldots, p_{n}$, $\Ical_{k}=I_{k}\cap\set{\absj{\omega},\ldots, t}$. On $\Omega_{coup}$, we have
\begin{align*}
\sum_{i=\absj{\omega}}^{t}Z_{i}&=\sum_{k=1-h_{1}}^{p_{n}}\sum_{i\in \Ical_{k}}Z_{i}=\sum_{k=1-h_{1}}^{p_{n}}\sum_{i\in \Ical_{k}}Z_{i}^{\prime}\enspace .
\end{align*}
Hence, for all $x>0$, $\nu\in(0,1)$, a union bound gives
\begin{multline}
\label{eq.conc.1}\P\set{\absj{\sum_{i=\absj{\omega}}^{t}Z_{i}}>x\cap\Omega_{coup}}
=\P\set{\absj{\sum_{k=0}^{k_n}\sum_{i\in \Ical_{1-h_{1}+2k}}Z_{i}^{\prime}}>\nu x\cap\Omega_{coup}}\\+\P\set{\absj{\sum_{k=0}^{\ell_{n}}\sum_{i\in \Ical_{2-h_{1}+2k}}Z_{i}^{\prime}}>(1-\nu)x\cap\Omega_{coup}}\\
\leq \P\set{\absj{\sum_{k=0}^{k_n}\sum_{i\in \Ical_{1-h_{1}+2k}}Z_{i}^{\prime}}>\nu x}+\P\set{\absj{\sum_{k=0}^{\ell_{n}}\sum_{i\in \Ical_{2-h_{1}+2k}}Z_{i}^{\prime}}>(1-\nu)x}\enspace .
\end{multline}
By construction,  $(\sum_{i\in \Ical_{1-h_{1}+2k}}Z_{i}^{\prime})_{k=0,\ldots,k_{n}}$ and $(\sum_{i\in \Ical_{2-h_{1}+2k}}Z_{i}^{\prime})_{k=0,\ldots,\ell_{n}}$ are independent and upper bounded by $(d_{n}+r_{n})/\sqrt{\mu(\omega)}$. Let 
\[
\sigma_{1}^{2}=\sum_{k=0}^{k_{n}}\var\paren{\sum_{i\in \Ical_{1-h_{1}+2k}}Z_{i}},\qquad \sigma_{2}^{2}=\sum_{k=0}^{\ell_{n}}\var\paren{\sum_{i\in \Ical_{2-h_{1}+2k}}Z_{i}}\enspace .
\]
From Lemma \ref{lem.covariance}, we have 
\[
\sigma_{1}^{2}\leq 2\paren{\Phi+\frac1{\lambda-1}}\mu(\omega)\sum_{k=0}^{k_{n}}\card\set{\Ical_{1-h_{1}+2k}}\enspace ,
\]
\[
\sigma_{2}^{2}\leq 2\paren{\Phi+\frac1{\lambda-1}}\mu(\omega)\sum_{k=0}^{\ell_{n}}\card\set{\Ical_{2-h_{1}+2k}}\enspace .
\]
Therefore, for $L=2\sqrt{\Phi+(\lambda-1)^{-1}}$,
\[
\;n_{1}=\sum_{k=0}^{k_{n}}\card\set{\Ical_{1-h_{1}+2k}}\;\mbox{and} \; n_{2}=\sum_{k=0}^{\ell_{n}}\card\set{\Ical_{2-h_{1}+2k}}\enspace ,
\] 
Benett's inequality (see Lemma \ref{lem.benett}) yields that for all $y>0$,
\begin{align*}
\P\set{\absj{\sum_{k=0}^{k_{n}}\sum_{i\in \Ical_{1-h_{1}+2k}}Z_{i}^{\prime}}>L\sqrt{n_{1}y}+\frac {(d_{n}+r_{n})y}{3\sqrt{\mu(\omega)}}}&\leq 2e^{-y},\\
\P\set{\absj{\sum_{k=0}^{\ell_{n}}\sum_{i\in \Ical_{2-h_{1}+2k}}Z_{i}^{\prime}}>L\sqrt{n_{2}y}+\frac {(d_{n}+r_{n})y}{3\sqrt{\mu(\omega)}}}&\leq 2e^{-y}\enspace .
\end{align*}
In \eqref{eq.conc.1}, we choose
\[
x=L\paren{\sqrt{n_{1}}+\sqrt{n_{2}}}\sqrt{y}+\frac {(d_{n}+r_{n})y}{3(\nu\wedge (1-\nu))\sqrt{\mu(\omega)}},\;\nu=\frac{\sqrt{n_{1}}}{\sqrt{n_{1}}+\sqrt{n_{2}}}
\]
We have $\sqrt{n_{1}}+\sqrt{n_{2}}\leq 2\sqrt{t-\absj{\omega}+1}$ and 
\[
\frac{n_{1}\wedge n_{2}}{n_{1}+n_{2}}\geq \frac{p_{n}-3+h_{1}}{2p_{n}-3+h_{1}}\geq \frac{p_{n}-1-2}{2(p_{n}-1)}\geq \frac12-\frac{1}{p_{n}-1}\geq \frac16\;.
\]
Hence, $(\nu\wedge (1-\nu))\geq \sqrt{(n_{1}\wedge n_{2})/(n_{1}+n_{2})}\geq 1/3$ and we have obtained that, for all $y>0$,
\[
\P\set{\absj{\sum_{i=\absj{\omega}}^{t}Z_{i}}>2L\sqrt{n-\absj{\omega}+1}\sqrt{y}+\frac {(d_{n}+r_{n})y}{\sqrt{\mu(\omega)}}\cap\Omega_{coup}}\leq 4e^{-y}\enspace .
\]
This result can be rewritten as~\eqref{eq.conc.bas}.
\subsection{A complement for slope heuristic in the mixing case}\label{section.slope.mix}
\begin{proposition}\label{prop.slope.mix}
Let $(X_n)_{n\in\Z}$ be a $\phi$-mixing process satisfying \eqref{cond.phi.mix}, \eqref{cond.ND} and \eqref{cond.geo.bet.mix}.
\begin{align*}
\Fcalst^{(n)}&=\set{\omega\in\Fcal,\telque \forall a\in A,\;(n-\absj{\omega}+1)\muh_{n}(\omega a)\in\set{0}\cup\left[(\ln n)^{4},+\infty\right[}\enspace.
\end{align*}
Let $\tau\prec \taust$ be two trees in $\Fcalst^{(n)}$. For any $k\geq d_{n}+1$, let us define 
\[
Z_{k}=\sum_{(\omega,a)\in \taust\times A}\frac{\un_{X_{k-|\omega|}^{k}=\omega a}-\mu(\omega a)}{n-|\omega|}\ln\paren{\frac{P_{\taust}(a|\omega)}{P_{\tau}(a|\omega_{\tau})}}\enspace .
\]
Let $\Omega_{coup}$ be the event defined in \eqref{def.Omega.coup}. For any $x>0$, $\epsilon\in (0,1)$, let $d_{n}=\log n$,
\begin{multline*}
\P\set{\absj{\sum_{k=d_{n}+1}^{n}Z_{k}}>\epsilon K_{\mu_{\taust}}(P_{\taust},P_{\tau})+L\frac{\epsilon^{-1}\Lst^{(\tau,\taust)}\vee (\ln n)^{2}}{(n-d_{n})}x}\\
\leq \P\set{\Omega_{coup}^{c}}+2e^{-x}\enspace .
\end{multline*}
\end{proposition}
\begin{remark}
$\sum_{k=d_{n}+1}^{n}Z_{k}$ is essentially equal to $L(P_{\taust})-L(P_{\tau})$. Proposition \ref{prop.slope.mix} and a union bound state then that, in the mixing case, the event $\Omega_{mis}$ defined in Theorem \ref{theo.opt.pen.afaire} holds if $\Fcalst^{(n)}$ is contained in a fixed set of trees with cardinality polynomial in $n$. 
\end{remark}

\begin{proof} Let us keep the notation of the proof of Theorem \ref{theo.conc.bas}.\\
$\mbox{If}\;\exists j\in\set{0,\ldots,\ell_{n}}\telque \set{k-d_{n}+1,\ldots, k}\subset J_{j}^{\prime}$, let
\[
Z^{\prime}_{k}=\sum_{(\omega,a)\in \taust\times A}\frac{\un_{(Y^{\prime})_{k-|\omega|}^{k}=\omega a}-\mu(\omega a)}{n-|\omega|}\ln\paren{\frac{P_{\taust}(a|\omega)}{P_{\tau}(a|\omega_{\tau})}}\enspace ,
\]
Otherwise, let
\[
Z^{\prime}_{k}= \sum_{(\omega,a)\in \taust\times A}\frac{\un_{Y_{k-|\omega|}^{k}=\omega a}-\mu(\omega a)}{n-|\omega|}\ln\paren{\frac{P_{\taust}(a|\omega)}{P_{\tau}(a|\omega_{\tau})}}\enspace .
\]
For any $j=0,\ldots,\kappa_{n}$, we denote by $\Ical_{j}$ the set of values of $k$ such that $\set{k-d_{n}+1,\ldots, k}\subset J_{j}$ and, for any $j=0,\ldots,\ell_{n}$, by $\Ical_{j}^{\prime}$ the set of $k$ such that $\set{k-d_{n}+1,\ldots, k}\subset J_{j}^{\prime}$. For any $x>0$ and $\nu\in (0,1)$, we have
\begin{multline*}
\P\set{\absj{\sum_{k=d_{n}+1}^{n}Z_{k}}>x}\leq \P\set{\Omega_{coup}^{c}}+\P\set{\absj{\sum_{j=0}^{\kappa_{n}}\sum_{k\in \Ical_{j}}Z_{k}}>\nu x\cap\Omega_{coup}}\\
+\P\set{\absj{\sum_{j=0}^{\ell_{n}}\sum_{k\in \Ical^{\prime}_{j}}Z_{k}}>(1-\nu)x\cap\Omega_{coup}}\\
\leq \P\set{\Omega_{coup}^{c}}+\P\set{\absj{\sum_{j=0}^{\kappa_{n}}\sum_{k\in \Ical_{j}}Z^{\prime}_{k}}>\nu x}+\P\set{\absj{\sum_{j=0}^{\ell_{n}}\sum_{k\in \Ical^{\prime}_{j}}Z^{\prime}_{k}}>(1-\nu)x}\enspace .
\end{multline*}
The random variables $(\sum_{k\in \Ical_{j}}Z^{\prime}_{k})_{j=0,\ldots,\kappa_{n}}$ and $(\sum_{k\in \Ical^{\prime}_{j}}Z^{\prime}_{k})_{j=0,\ldots,\ell_{n}}$ are independent by construction. Therefore, Benett's inequality yields, for any $x>0$,
\begin{equation}\label{eq.Benn.ineq.bas}
\P\set{\absj{\sum_{j=0}^{\kappa_{n}}\sum_{k\in \Ical_{j}}Z^{\prime}_{k}}>\sqrt{2\sigma_{1}^{2}x}+\frac{bx}3}\leq 2e^{-x}\enspace .
\end{equation}
In the previous inequality,
\[
\sigma_{1}^{2}\geq \sum_{j=1}^{\kappa_{n}}\var\paren{\sum_{k\in \Ical_{j}}Z_{k}} \hbox{ and }\qquad b\geq \max_{j=1,\ldots,\kappa_{n}}\norm{\sum_{k\in \Ical_{j}}Z^{\prime}_{k}}_{\infty}\enspace .
\]
The typicality property implies that, for $n$ large enough,
\begin{align}
\notag \norm{\sum_{k\in \Ical_{j}}Z^{\prime}_{k}}_{\infty}&\leq \max_{(\omega,a)\in\taust\times A\telque P(a|\omega)\neq 0}\frac{r_{n}+d_{n}}{n-d_{n}}\ln\paren{\frac1{P(a|\omega)}}\\
\label{eq.cont.norm.sup}&\leq \frac{(r_{n}+d_{n})\ln (2n)}{n-d_{n}}\enspace .
\end{align}
Moreover, for any $u\geq d_{n}+1$, by stationarity of $X_{1}^{n}$,
\begin{multline*}
\var\paren{\sum_{k\in \Ical_{j}}Z_{k}}=\sum_{k\in \Ical_{j}}\var(Z_{k})+\sum_{k\neq k^{\prime}\in \Ical_{j}}\cov(Z_{k},Z_{k^{\prime}})\\
=\card\set{\Ical_{j}}\var(Z_{u})+2\sum_{k=u+1}^{u+r_{n}}(r_{n}+u-k+1)\cov\paren{Z_{u},Z_{k}}\enspace .
\end{multline*}
We have
\begin{align*}
\var&(Z_{u})\leq \E\paren{\paren{\sum_{(\omega,a)\in \taust\times A}\frac{\un_{X_{k-|\omega|}^{k}=\omega a}}{n-|\omega|}\ln\paren{\frac{P_{\taust}(a|\omega)}{P_{\tau}(a|\omega_{\tau})}}}^{2}}\\
&\leq \sum_{(\omega,a)\in \taust\times A}\paren{\ln\paren{\frac{P_{\taust}(a|\omega)}{P_{\tau}(a|\omega_{\tau})}}}^{2}\frac{\mu(\omega a)}{(n-d_{n})^{2}}\\
&\leq \Lst^{(\tau,\taust)}\sum_{\omega\in \taust}\frac{\mu(\omega)}{(n-d_{n})^{2}}\sum_{a\in A}P_{\taust}(a|\omega)\wedge P_{\tau}(a|\omega_{\tau}) \paren{\ln\paren{\frac{P_{\taust}(a|\omega)}{P_{\tau}(a|\omega_{\tau})}}}^{2}.
\end{align*}
Lemma \ref{lem:BarSheu} gives 
\begin{equation}\label{eq.var.int}
\var(Z_{u})\leq \frac{\Lst^{(\tau,\taust)}}{(n-d_{n})^{2}}K_{\mu_{\taust}}(P_{\taust},P_{\tau})\enspace .
\end{equation}
In addition, using Lemma \ref{lem.covariance}, we get, for 
\[
m_{+}(u,k)=\max_{\omega\in \taust}\set{\sqrt{\phi_{k-u-\absj{\omega}+1}}\un_{k-u-\absj{\omega}+1\geq 0}+\lambda^{k-u}\un_{k-u-\absj{\omega}+1< 0}}\enspace ,
\]
\begin{align}
\label{cov.ineq.1}(n-&d_{n})^{2}\cov\paren{Z_{u},Z_{k}}\\
&\leq
\notag\sum_{\set{(\omega,a),(\omega^{\prime},a^{\prime})}\in (\taust\times A)^{2}} \absj{\ln\paren{\frac{P_{\taust}(a|\omega)}{P_{\tau}(a|\omega_{\tau})}}}
 \absj{\ln\paren{\frac{P_{\taust}(a^{\prime}|\omega^{\prime})}{P_{\tau}(a^{\prime}|\omega^{\prime}_{\tau})}}} \\
&\quad\quad\quad\times \absj{\cov\paren{\un_{X_{u-|\omega^{\prime}|}^{u}=\omega^{\prime}a},\un_{X_{k-\absj{\omega }}^{k}=\omega a}}}\\
\notag&\leq m_{+}(u,k)\paren{\sum_{(\omega,a)\in \taust\times A}\sqrt{\mu(\omega a)}\absj{\ln\paren{\frac{P_{\taust}(a|\omega)}{P_{\tau}(a|\omega_{\tau})}}}}^{2}\enspace .
\end{align}
The Cauchy-Schwarz inequality yields
\begin{align*}
&\left(\sum_{(\omega,a)\in \taust\times A}\right.\left.\sqrt{\mu(\omega a)}\absj{\ln\paren{\frac{P_{\taust}(a|\omega)}{P_{\tau}(a|\omega_{\tau})}}}\right)^{2}\\
&\leq N(\taust)|A|\sum_{\omega\in \taust}\mu(\omega)\sum_{a\in A}P_{\taust}(a|\omega)\paren{\ln\paren{\frac{P_{\taust}(a|\omega)}{P_{\tau}(a|\omega_{\tau})}}}^{2}\\
&\leq \Lst^{(\tau,\taust)}N(\taust)|A|\sum_{\omega\in \taust}\mu(\omega)\sum_{a\in A}P_{\taust}(a|\omega)\wedge P_{\tau}(a|\omega_{\tau}) \paren{\ln\paren{\frac{P_{\taust}(a|\omega)}{P_{\tau}(a|\omega_{\tau})}}}^{2}\enspace .
\end{align*}
Using Lemma \ref{lem:BarSheu},
\[
\paren{\sum_{(\omega,a)\in \taust\times A}\sqrt{\mu(\omega a)}\absj{\ln\paren{\frac{P_{\taust}(a|\omega)}{P_{\tau}(a|\omega_{\tau})}}}}^{2}\leq  \Lst^{(\tau,\taust)}N(\taust)|A|K_{\mu_{\taust}}(P_{\taust},P_{\tau})\enspace .
\]
Plugging this inequality in \eqref{cov.ineq.1} gives
\[
\cov\paren{Z_{u},Z_{k}}\leq m_{+}(u,k)\Lst^{(\tau,\taust)}\frac{N(\taust)|A|}{(n-d_{n})^{2}}K_{\mu_{\taust}}(P_{\taust},P_{\tau})\enspace .
\]
As $|\omega|\leq d_{n}$, we have, under \eqref{cond.geo.bet.mix}, $m_{+}(u,k)\leq \paren{L_{mix}e^{\gamma_{mix}d_{n}}\vee 1}(\lambda\vee e^{-\gamma_{mix}})^{k-u}$.

We always also have the basic inequality
\[
\cov\paren{Z_{u},Z_{k}}\leq \var\paren{Z_{u}}\leq \frac{\Lst^{(\tau,\taust)}}{(n-d_{n})^{2}}K_{\mu_{\taust}}(P_{\taust},P_{\tau})\enspace .
\]
Therefore
\begin{align*}
&\var\paren{\sum_{k\in \Ical_{j}}Z_{k}}=\card\set{\Ical_{j}}\var(Z_{u})\\
&\quad\quad\quad\quad\quad\quad+2\sum_{k=u+1}^{u+r_{n}}(r_{n}+u-k+1)\cov\paren{Z_{u},Z_{k}}\\
&\leq 2\frac{\Lst^{(\tau,\taust)(r_{n}+d_{n})}}{(n-d_{n})^{2}}K_{\mu_{\taust}}(P_{\taust},P_{\tau})\\
&\quad\quad\quad\quad\quad\quad\times\sum_{k=0}^{\infty}\paren{1\wedge \set{N(\taust)|A|\paren{L_{mix}e^{\gamma_{mix}d_{n}}\vee 1}(\lambda\vee e^{-\gamma_{mix}})^{k}}}\enspace .
\end{align*}
As $N(\taust)\leq n$ and $e^{\gamma_{mix}d_{n}}\leq n^{\gamma_{mix}/\ln(|A|)}$, we can cut the separate between the $k\leq \ln\paren{ n^{1+\gamma_{mix}/\ln(|A|)}} / \ln (\lambda\vee e^{-\gamma_{mix}})$ and the $k>\ln\paren{ n^{1+\gamma_{mix}/\ln(|A|)}} /\ln (\lambda\vee e^{-\gamma_{mix}})$, and we obtain that there exists a constant $L\egaldef L(L_{mix}, \gamma_{mix}, |A|, \lambda)$ such that
\[
\var\paren{\sum_{k\in \Ical_{j}}Z_{k}}\leq Lr_{n}\frac{\Lst^{(\tau,\taust)}}{(n-d_{n})^{2}}K_{\mu_{\taust}}(P_{\taust},P_{\tau})
\]
Plugging this inequality and \eqref{eq.cont.norm.sup} in \eqref{eq.Benn.ineq.bas} gives
\begin{multline*}
\P\set{\absj{\sum_{j=0}^{\kappa_{n}}\sum_{k\in \Ical_{j}}Z^{\prime}_{k}}>\sqrt{L\kappa_{n}r_{n}\frac{\Lst^{(\tau,\taust)}}{(n-d_{n})^{2}}K_{\mu_{\taust}}(P_{\taust},P_{\tau})x}+\frac{r_{n}\ln (2n)}{n-d_{n}}\frac{x}3}\\\leq 2e^{-x}\enspace .
\end{multline*}
Using the basic inequality $2ab\leq \epsilon a^{2}+\epsilon^{-1}b^{2}$, we finally get
\[
\P\set{\absj{\sum_{j=0}^{\kappa_{n}}\sum_{k\in \Ical_{j}}Z^{\prime}_{k}}>\frac{\epsilon}2K_{\mu_{\taust}}(P_{\taust},P_{\tau})+L\frac{\epsilon^{-1}\Lst^{(\tau,\taust)}\vee (\ln n)^{2}}{(n-d_{n})}x}\leq 2e^{-x}\enspace .
\]
The same control holds for $\absj{\sum_{j=0}^{\ell_{n}}\sum_{k\in \Ical^{\prime}_{j}}Z^{\prime}_{k}}$, hence
\begin{multline*}
\P\set{\absj{\sum_{k=d_{n}+1}^{n}Z_{k}}>\epsilon K_{\mu_{\taust}}(P_{\taust},P_{\tau})+L\frac{\epsilon^{-1}\Lst^{(\tau,\taust)}\vee (\ln n)^{2}}{(n-d_{n})}x}\\
\leq \P\set{\Omega_{coup}^{c}}+2e^{-x}\enspace .
\end{multline*}
\end{proof}

\section{Links with continuity rates}\label{sec.cont.rates}
\subsection{Control of the bias with the continuity rates}
An important tool in the theory of chains of infinite order is the continuity rates defined, $\forall (a_{-\infty}^{-1},a,k)\in A^{-\N}\times A\times \N^{*},$ by
\begin{align*}
\epsilon_{k}(a_{-\infty}^{-1}a)\egaldef\sup_{(b_{-\infty}^{-k-1},c_{-\infty}^{-k-1})\in (A^{-\N})^{2}}\absj{P(a|b_{-\infty}^{-k-1}a_{-k}^{-1})-P(a|c_{-\infty}^{-k-1}a_{-k}^{-1})}\enspace .
\end{align*}
Let us remark that, for all $(a_{-\infty}^{-1},a,k)\in A^{-\N}\times A\times \N^{*}$, $\epsilon_{k}(a_{-\infty}^{-1}a)$ only depends on $(a_{-k}^{-1},a)$, therefore, we will also use the following notation
\[
\forall (a_{-\infty}^{-1},\omega,a)\in A^{-\N}\times A^{*}\times A,\qquad\epsilon(\omega,a)\egaldef\epsilon_{\absj{\omega}}(a_{-\infty}^{-1}\omega,a)\enspace .
\]
These continuity rates can be used to upper bound the bias term of the risk. In order to see this, we introduce the following definition
\[
\forall \tau\in \Tcal,\;\norm{\epsilon}_{\tau}^{2}\egaldef\int_{A^{-\N}}d\mu(\omega)\sum_{a\in A}\epsilon_{\absj{\omega_{\tau}}}(\omega a)^{2}\enspace .
\]

\begin{proposition}\label{prop.control.bias}
Let $\tau$ be a finite context tree, let $\eta\in(0,e^{-1})$ and let 
\[
\Omega_{\eta}=\set{\omega\in A^{-\N} \telque \exists a\in A,\;\olP_{\tau}(a|\omega)< \eta}\enspace .
\]
Then,
\begin{equation}\label{eq.control.bias}
K_{\mu}\paren{P,\olP_{\tau}}\leq \frac{\norm{\epsilon}_{\tau}^{2}}{\eta}+|A|\eta\ln\paren{\frac1{\eta}}\mu\paren{\Omega_{\eta}}\enspace.
\end{equation}
\end{proposition}

Proposition \eqref{prop.control.bias} can be used under the following assumption
\begin{equation}\label{Hyp.Gibbs}
\hyptag{GC}\exists \Kst>0\telque \forall (\omega,a)\in A^{-\N}\times A,\qquad P(a|\omega)\geq \frac1{\Kst}\enspace .
\end{equation}

In that case, $\Omega_{\eta}=\emptyset$ for all $\eta<\Kst^{-1}$, hence \eqref{eq.control.bias} yields 
\[
K_{\mu}\paren{P,\olP_{\tau}}\leq \Kst\norm{\epsilon}_{\tau}^{2}\enspace .
\]
If, on the other hand, $\mu(\Omega_{\eta})>0$ for all $\eta>0$, we can choose \[\eta=\norm{\epsilon}_{\tau}\paren{\mu\paren{\Omega_{\norm{\epsilon}_{\tau}}}}^{-1/2}\] and get an absolute constant $C$ such that, for all $r\in(0,1)$, 
\[
K_{\mu}\paren{P,\olP_{\tau}}\leq \frac{C}{r}\paren{\norm{\epsilon}_{\tau}\sqrt{\mu\paren{\Omega_{\norm{\epsilon}_{\tau}}}}}^{1-r}\enspace .
\]

\begin{proof}
By definition, for all $(\omega,a)\in A^{-\N}\times A$,
\begin{equation}\label{Proof.control.Bias.1}
\absj{P(a|\omega)-\olP_{\tau}(a|\omega)}\leq \epsilon_{\absj{\omega_{\tau}}}(\omega a)\enspace .
\end{equation}
In addition, we have
\begin{equation}\label{Proof.control.Bias.2}
K_{\mu}\paren{P,\olP_{\tau}}=\int_{A^{-\N}\times A}d\mu(\omega)P(a|\omega)\ln\paren{\frac{P(a|\omega)}{\olP_{\tau}(a|\omega)}}\enspace .
\end{equation}
Let $\eta\in(0,e^{-1})$ and let 
\begin{align*}
\Omega_{1,\eta}&=\set{(\omega,a)\in A^{-\N}\times A \telque \olP_{\tau}(a|\omega)\geq \eta}\enspace ,\\
\Omega_{2,\eta}&=\set{(\omega,a)\notin \Omega_{1,\eta} \telque \olP_{\tau}(a|\omega)<P(a|\omega)}\enspace ,\\
\Omega_{3,\eta}&=\set{(\omega,a)\in \tau\times A \telque P_{\tau}(a|\omega)< \eta}\enspace .
\end{align*}
From \eqref{Proof.control.Bias.2}, we have
\begin{align*}
K_{\mu}\paren{P,\olP_{\tau}}&=\paren{\int_{\Omega_{1,\eta}}+\int_{\Omega_{1,\eta}^{c}}}d\mu(\omega)P(a|\omega)\ln\paren{\frac{P(a|\omega)}{\olP_{\tau}(a|\omega)}}\\
&\leq \paren{\int_{\Omega_{1,\eta}}+\int_{\Omega_{2,\eta}}}d\mu(\omega)P(a|\omega)\ln\paren{\frac{P(a|\omega)}{\olP_{\tau}(a|\omega)}}\\
&\leq \int_{\Omega_{1,\eta}}d\mu(\omega)P(a|\omega)\ln\paren{\frac{P(a|\omega)}{\olP_{\tau}(a|\omega)}}\\
&\quad\quad\quad\quad+\int_{\Omega_{2,\eta}} d\mu(\omega)P(a|\omega)\ln\paren{\frac{1}{\olP_{\tau}(a|\omega)}}\\
&\leq \int_{\Omega_{1,\eta}}d\mu(\omega)P(a|\omega)\ln\paren{\frac{P(a|\omega)}{\olP_{\tau}(a|\omega)}}\\
&\quad\quad\quad\quad+\sum_{(\omega,a)\in\Omega_{3,\eta}}\mu_{\tau}(\omega)P_{\tau}(a|\omega)\ln\paren{\frac{1}{P_{\tau}(a|\omega)}}\enspace .
\end{align*}
We use the bound
\[
\forall x\leq \eta,\;x\ln\paren{\frac1x}\leq \eta\ln\paren{\frac1{\eta}}\enspace .
\]
We obtain
\[
\sum_{(\omega,a)\in\Omega_{3,\eta}}\mu_{\tau}(\omega)P_{\tau}(a|\omega)\ln\paren{\frac{1}{P_{\tau}(a|\omega)}}\leq |A|\eta\ln\paren{\frac1{\eta}}\mu\set{\Omega_{\eta}}\enspace .
\]
In addition, since $\Omega_{1,\eta}$ does not depend on the pasts before $\tau$, $\mu(\Omega_{1,\eta})=\mu_{\tau}(\Omega_{1,\eta})$ and, using that $\forall x>0,\;\ln(x)\leq x-1$, we obtain
\begin{align*}
\int_{\Omega_{1,\eta}}&d\mu(\omega)P(a|\omega)\ln\paren{\frac{P(a|\omega)}{\olP_{\tau}(a|\omega)}}\\
&=\int_{\Omega_{1,\eta}}d\mu(\omega)\paren{P(a|\omega)-\olP_{\tau}(a|\omega)+\olP_{\tau}(a|\omega)}\paren{\frac{P(a|\omega)-\olP_{\tau}(a|\omega)}{\olP_{\tau}(a|\omega)}}\\
&\leq \frac1{\eta}\int_{A^{\N}}d\mu(\omega)\sum_{a\in A}\paren{P(a|\omega)-\olP_{\tau}(a|\omega)}^{2}+\mu(\Omega_{1,\eta})-\mu_{\tau}(\Omega_{1,\eta})\\
&\leq\frac1{\eta}\int_{A^{\N}}d\mu(\omega)\sum_{a\in A}\epsilon_{\absj{\omega_{\tau}}}(\omega a)^{2}=\frac{\norm{\epsilon}_{\tau}^{2}}{\eta}\enspace .
\end{align*}
\end{proof}

\subsection{Mixing properties and continuity rates}\label{section.continuity.mixing}
$\phi$-mixing conditions can also be deduced from continuity. In order to see that, let us recall the following equivalent definition of $\phi$-mixing coefficient (see \cite{Br02} prop 3.22)
\[
\phi(k)=\sup_{s\in \N}\sup_{E\in A^{s}}\norm{P(X_{k}^{k+s}\in E|\sigma(X_{-\infty}^{0}))-P(E)}_{\infty}\enspace .
\]
Let us introduce the following assumptions.
\begin{equation}\label{cond.cont.expo}
\hyptag{EC}\exists (C,\alpha)\in (\R_{+}^{*})>0\telque \forall \ell>0, \; 1-\inf_{\omega^{\prime}\in A^{\ell}}\sum_{a\in A}\inf_{\omega\in A^{-\N}}P(a|\omega\omega^{\prime})\leq Ce^{-\alpha\ell}
\end{equation}
\begin{equation}\label{cond.regeneration}
\hyptag{RC}\exists p_{\min}>0 \telque \forall (a,\omega)\in A\times A^{-\N},\; P(a|\omega)>0\enspace .
\end{equation}
From Theorem 4.1 and Corollary 4.1 in \cite{CFF02}, under assumptions \eqref{cond.cont.expo} and \eqref{cond.regeneration}, there exists $C_{i}>0$, $\alpha_{i}>0$ such that
\begin{multline*}
\sup_{s\in \N}\sup_{E\in A^{s}}\norm{P(X_{k+1}^{k+s}\in E|\sigma(X_{-\infty}^{0}))-P(E)}_{\infty}\\\leq 2C_{i}\sum_{j=0}^{\infty}e^{-\alpha_{i}(k+j)}\leq \frac{2C_{i}}{1-e^{-\alpha_{i}}}e^{-\alpha_{i}k}\enspace .
\end{multline*}
Let then $\norm{\epsilon}_{k,\infty}=\sup_{(\omega,a)\in A^{-\N}\times A}\epsilon_{k}(\omega, a)$. It is clear that 
\[
1-\inf_{\omega^{\prime}\in A^{k}}\sum_{a\in A}\inf_{\omega\in A^{-\N}}P(a|\omega\omega^{\prime})\leq \norm{\epsilon}_{k,\infty}\enspace .
\] 
Therefore, we have proved the following proposition.
\begin{proposition}
Every stationary ergodic process satisfying \eqref{cond.regeneration} and such that $\norm{\epsilon}_{k,\infty}$ decreases exponentially satisfies Assumption \eqref{cond.geo.bet.mix}.
\end{proposition}

\section{Technical tools}
In the main proofs, we used the following lemmas.
\subsection{Decomposition of the risk}
\begin{lemma}\label{lem.dec.oracle}
For all $\tau\in \Fcal$, 
\begin{align*}
\int_{A^{-\N}\times A}d\mu(\omega a)\ln\paren{P(a|\omega)}&+\sum_{(\omega,a)\in \tau\times A}\muh_{n}(\omega a)\ln\paren{\frac1{\Ph_{\tau}(a|\omega)}}\\
&=K_{\mu}(P,\olP_{\tau})-K_{\muh}(\Ph_{\tau},P_{\tau})+L(P_{\tau})\enspace ,
\end{align*}
where
\[
L(P_{\tau})\egaldef \sum_{(\omega,a)\in \tau\times A}(\muh_{n}(\omega a)-\mu(\omega a))\ln\paren{\frac1{P_{\tau}(a|\omega)}}\enspace .
\]
\end{lemma}
\begin{proof}
\begin{align*}
\int_{A^{-\N}\times A}&d\mu(\omega a)\ln\paren{P(a|\omega)}+\sum_{(\omega,a)\in \tau\times A}\muh_{n}(\omega a)\ln\paren{\frac1{\Ph_{\tau}(a|\omega)}}\\
&=K_{\mu}(P,\olP_{\tau})-\sum_{(\omega,a)\in \tau\times A}\mu(\omega a)\ln\paren{\frac1{P_{\tau}(a|\omega)}}\\
&\quad\quad\quad\quad+\sum_{(\omega,a)\in \tau\times A}\muh_{n}(\omega a)\ln\paren{\frac1{\Ph_{\tau}(a|\omega)}}\\
&=K_{\mu}(P,\olP_{\tau})+L(P_{\tau})+\sum_{(\omega,a)\in \tau\times A}\muh_{n}(\omega a)\ln\paren{\frac{P_{\tau}(a|\omega)}{\Ph_{\tau}(a|\omega)}}\\
&=K_{\mu}(P,\olP_{\tau})+L(P_{\tau})-K_{\muh}(\Ph_{\tau},P_{\tau})\enspace .
\end{align*}
\end{proof}

\subsection{Control of $L(P_{\tau})-L(P_{\tau^{\prime}})$}
\begin{lemma}\label{lem.termes.croises.1}
For all $(\tau,\tau^{\prime})\in \Fcal^{2}$, let $\Tcal(\tau,\tau^{\prime})$ be the unique tree satisfying the following conditions.
\begin{enumerate}
\item $\tau\prec\Tcal(\tau,\tau^{\prime})$ and $\tau^{\prime}\prec\Tcal(\tau,\tau^{\prime})$.
\item $\Tcal(\tau,\tau^{\prime})\subset \tau\cup\tau^{\prime}$.
\end{enumerate}
Then,
\[
L(P_{\tau})-L(P_{\tau^{\prime}})=\sum_{(\omega,a)\in \Tcal(\tau,\tau^{\prime})\times A}(\muh_{n}(\omega a)-\mu(\omega a))\ln\paren{\frac{P_{\tau^{\prime}}(a|\omega_{\tau^{\prime}})}{P_{\tau}(a|\omega_{\tau})}}\enspace .
\]
\end{lemma}
\begin{proof}
The result follows from the following remark.
Let $\tau\in \Fcal$ and let $\overline{\tau}$ be any element of $\Fcal$ such that $\tau\prec\overline{\tau}$. As $\mu$ and $\muh$ are probability measures, we have
\[
L(P_{\tau})=\sum_{(\omega,a)\in \overline{\tau}\times A}(\muh_{n}(\omega a)-\mu(\omega a))\ln\paren{\frac{1}{P_{\tau}(a|\omega_{\tau})}}\enspace .
\]
\end{proof}

\begin{lemma}\label{lem.termes.croises.2}
Let $(\tau,\tau^{\prime})\in \Fcal^{2}$ and let $\Tcal(\tau,\tau^{\prime})$ be the associated tree defined in Lemma \ref{lem.termes.croises.1}. Let 
\begin{align*}
\eta&=\max_{\omega\in \Tcal(\tau,\tau^{\prime}),\mu(\omega)\neq 0}\set{\frac{\absj{\muh_{n}(\omega)-\mu(\omega )}}{\mu(\omega )}}\\
\Lst^{(\tau,\tau^{\prime})}&=\max_{(\omega,a)\in \Tcal(\tau,\tau^{\prime})\times A,\;P_{\tau}(a|\omega_{\tau})\wedge P_{\tau^{\prime}}(a|\omega_{\tau^{\prime}})\neq 0}\set{\frac{P_{\tau}(a|\omega_{\tau})\vee P_{\tau^{\prime}}(a|\omega_{\tau^{\prime}})}{P_{\tau}(a|\omega_{\tau})\wedge P_{\tau^{\prime}}(a|\omega_{\tau^{\prime}})}}\enspace .
\end{align*}
Then, for all $\epsilon>0$
\begin{align*}
L(P_{\tau})-L(P_{\tau^{\prime}})\leq& (\epsilon+O(\eta))\paren{K_{\mu}(P,\olP_{\tau})+K_{\mu}(P,\olP_{\tau^{\prime}})}\\
&+(1+O(\eta))\frac{\Lst^{(\tau,\tau^{\prime})}}{\epsilon}\paren{K_{\mu_{\tau}}(P_{\tau},\Ph_{\tau})+K_{\mu_{\tau^{\prime}}}(P_{\tau^{\prime}},\Ph_{\tau^{\prime}})}\enspace .
\end{align*}
\end{lemma}

\begin{proof}
From Lemma \ref{lem.termes.croises.1}, we have
\begin{align}
\notag L(P_{\tau})&-L(P_{\tau^{\prime}})= \\\notag&\sum_{\omega\in \Tcal(\tau,\tau^{\prime})}(\muh_{n-1}(\omega)-\mu(\omega))\sum_{a\in A}P_{\Tcal(\tau,\tau^{\prime})}(a|\omega)\ln\paren{\frac{P_{\Tcal(\tau,\tau^{\prime})}(a|\omega)}{P_{\tau}(a|\omega_{\tau})}}\\
\notag &+ \sum_{\omega\in \Tcal(\tau,\tau^{\prime})}(\mu(\omega)-\muh_{n-1}(\omega))\sum_{a\in A}P_{\Tcal(\tau,\tau^{\prime})}(a|\omega)\ln\paren{\frac{P_{\Tcal(\tau,\tau^{\prime})}(a|\omega)}{P_{\tau^{\prime}}(a|\omega_{\tau^{\prime}})}}\\
\label{eq.dec.L.1}&+ \sum_{\omega\in \Tcal(\tau,\tau^{\prime})}\muh_{n-1}(\omega)\sum_{a\in A}(P_{\Tcal(\tau,\tau^{\prime})}(a|\omega)-\Ph(a|\omega))\ln\paren{\frac{P_{\tau^{\prime}}(a|\omega_{\tau^{\prime}})}{P_{\tau}(a|\omega_{\tau})}}\enspace.
\end{align}
We have, for $\tau^{*}=\tau$ or $\tau^{\prime}$ 
\begin{align}
&\notag \absj{\sum_{\omega\in \Tcal(\tau,\tau^{\prime})}(\muh_{n-1}(\omega)-\mu(\omega))\sum_{a\in A}P_{\Tcal(\tau,\tau^{\prime})}(a|\omega)\ln\paren{\frac{P_{\Tcal(\tau,\tau^{\prime})}(a|\omega)}{P_{\tau}(a|\omega_{\tau^{*}})}}}\leq\\
&\notag  \eta\sum_{\omega\in \Tcal(\tau,\tau^{\prime})}\mu(\omega)\sum_{a\in A}P_{\Tcal(\tau,\tau^{\prime})}(a|\omega)\ln\paren{\frac{P_{\Tcal(\tau,\tau^{\prime})}(a|\omega)}{P_{\tau}(a|\omega_{\tau^{*}})}}\\
&\label{eq.dec.2}=\eta K_{\mu_{\Tcal(\tau,\tau^{\prime})}}(P_{\Tcal(\tau,\tau^{\prime})},P_{\tau^{*}})
\end{align}
Hence, in \eqref{eq.dec.L.1}, we have
\begin{align*}
&\sum_{\omega\in \Tcal(\tau,\tau^{\prime})}(\muh_{n-1}(\omega)-\mu(\omega))\sum_{a\in A}P_{\Tcal(\tau,\tau^{\prime})}(a|\omega)\ln\paren{\frac{P_{\Tcal(\tau,\tau^{\prime})}(a|\omega)}{P_{\tau}(a|\omega_{\tau})}}\\
&+ \sum_{\omega\in \Tcal(\tau,\tau^{\prime})}(\mu(\omega)-\muh_{n-1}(\omega))\sum_{a\in A}P_{\Tcal(\tau,\tau^{\prime})}(a|\omega)\ln\paren{\frac{P_{\Tcal(\tau,\tau^{\prime})}(a|\omega)}{P_{\tau^{\prime}}(a|\omega_{\tau^{\prime}})}}\\
&\leq \eta\paren{K_{\mu}(P,\olP_{\tau})+K_{\mu}(P,\olP_{\tau^{\prime}})}\enspace .
\end{align*}
Moreover, 
\begin{align}
\notag&\sum_{\omega\in \Tcal(\tau,\tau^{\prime})}\muh_{n-1}(\omega)\sum_{a\in A}(P_{\Tcal(\tau,\tau^{\prime})}(a|\omega)-\Ph(a|\omega))\ln\paren{\frac{P_{\tau^{\prime}}(a|\omega_{\tau^{\prime}})}{P_{\tau}(a|\omega_{\tau})}}\\
\label{eq.dec.3}&=\sum_{\omega\in \Tcal(\tau,\tau^{\prime})}\muh_{n-1}(\omega)\sum_{a\in A}(P_{\Tcal(\tau,\tau^{\prime})}(a|\omega)-\Ph(a|\omega))\\
&\quad\quad\quad\times\paren{\ln\paren{\frac{P_{\Tcal(\tau,\tau^{\prime})}(a|\omega)}{P_{\tau}(a|\omega_{\tau})}}-\ln\paren{\frac{P_{\Tcal(\tau,\tau^{\prime})}(a|\omega)}{P_{\tau^{\prime}}(a|\omega_{\tau^{\prime}})}}}
\end{align}
By Cauchy-Schwarz inequality, we have, for $\tau^{*}=\tau$ or $\tau^{\prime}$,
\begin{align}
\notag&\sum_{a\in A}(P_{\Tcal(\tau,\tau^{\prime})}(a|\omega)-\Ph(a|\omega))\ln\paren{\frac{P_{\Tcal(\tau,\tau^{\prime})}(a|\omega)}{P_{\tau^{*}}(a|\omega_{\tau^{*}})}} \leq\\
&\label{eq.dec.4}\sqrt{\sum_{a\in A}\frac{(P_{\Tcal(\tau,\tau^{\prime})}(a|\omega)-\Ph(a|\omega))^{2}}{P_{\Tcal(\tau,\tau^{\prime})}(a|\omega)}\sum_{a\in A}P_{\Tcal(\tau,\tau^{\prime})}(a|\omega)\paren{\ln\paren{\frac{P_{\Tcal(\tau,\tau^{\prime})}(a|\omega)}{P_{\tau^{*}}(a|\omega_{\tau^{*}})}}}^{2}}\enspace .
\end{align}
%\begin{align}
%\notag&\leq (1+\eta)\sum_{\omega\in \Tcal(\tau,\tau^{\prime})}\mu(\omega)\sqrt{\sum_{a\in A}\frac{(P_{\Tcal(\tau,\tau^{\prime})}(a|\omega)-\Ph(a|\omega))^{2}}{P_{\Tcal(\tau,\tau^{\prime})}(a|\omega)}\sum_{a\in A}P_{\Tcal(\tau,\tau^{\prime})}(a|\omega)\paren{\ln\paren{\frac{P_{\Tcal(\tau,\tau^{\prime})}(a|\omega)}{P_{\tau}(a|\omega_{\tau})}}}^{2}}\\
%&+(1+\eta)\sum_{\omega\in \Tcal(\tau,\tau^{\prime})}\mu(\omega)\sqrt{\sum_{a\in A}\frac{(P_{\Tcal(\tau,\tau^{\prime})}(a|\omega)-\Ph(a|\omega))^{2}}{P_{\Tcal(\tau,\tau^{\prime})}(a|\omega)}\sum_{a\in A}P_{\Tcal(\tau,\tau^{\prime})}(a|\omega)\paren{\ln\paren{\frac{P_{\Tcal(\tau,\tau^{\prime})}(a|\omega)}{P_{\tau^{\prime}}(a|\omega_{\tau^{\prime}})}}}^{2}}\enspace .
%\end{align}
%
Since $\Tcal(\tau,\tau^{\prime})\subset \tau\cup\tau^{\prime}$, we have
\begin{multline*}
\sum_{a\in A}\frac{(P_{\Tcal(\tau,\tau^{\prime})}(a|\omega)-\Ph(a|\omega))^{2}}{P_{\Tcal(\tau,\tau^{\prime})}(a|\omega)}\\\leq \sum_{a\in A}\frac{(P_{\tau}(a|\omega)-\Ph(a|\omega))^{2}}{P_{\tau}(a|\omega)}+\sum_{a\in A}\frac{(P_{\tau^{\prime}}(a|\omega)-\Ph(a|\omega))^{2}}{P_{\tau^{\prime}}(a|\omega)}\enspace .
\end{multline*}
From \eqref{eq.cons.typ}, \eqref{eq.cont.ln} and \eqref{eq.dec.Kull.1} in the proof of Lemma \ref{lem.cons.typ.gal}, we have, for $\tau^{*}=\tau$ or $\tau^{\prime}$
\[
\sum_{a\in A}\frac{(P_{\tau^{*}}(a|\omega)-\Ph(a|\omega))^{2}}{P_{\tau^{*}}(a|\omega)}\leq 2(1+O(\eta))\sum_{a\in A}P_{\tau^{*}}(a|\omega)\ln\paren{\frac{P_{\tau^{*}}(a|\omega)}{\Ph_{\tau^{*}}(a|\omega)}}\enspace .
\]
In addition, since $\Tcal(\tau,\tau^{\prime})\subset \tau\cup\tau^{\prime}$, for $\tau^{*}=\tau$ or $\tau^{\prime}$, we have
\begin{multline*}
\sum_{a\in A}P_{\Tcal(\tau,\tau^{\prime})}(a|\omega)\paren{\ln\paren{\frac{P_{\Tcal(\tau,\tau^{\prime})}(a|\omega)}{P_{\tau^{*}}(a|\omega_{\tau^{*}})}}}^{2}\\
\leq \Lst^{(\tau,\tau^{\prime})}\sum_{a\in A}\paren{P_{\Tcal(\tau,\tau^{\prime})}(a|\omega)\wedge P_{\tau^{*}}(a|\omega_{\tau^{*}})}\paren{\ln\paren{\frac{P_{\Tcal(\tau,\tau^{\prime})}(a|\omega)}{P_{\tau^{*}}(a|\omega_{\tau^{*}})}}}^{2}\enspace .
\end{multline*}
From Lemma \ref{lem:BarSheu}, we obtain 
\begin{align*}
\sum_{a\in A}P_{\Tcal(\tau,\tau^{\prime})}(a|\omega)&\paren{\ln\paren{\frac{P_{\Tcal(\tau,\tau^{\prime})}(a|\omega)}{P_{\tau^{*}}(a|\omega_{\tau^{*}})}}}^{2}\\
&\leq 2\Lst^{(\tau,\tau^{\prime})}\sum_{a\in A}P_{\Tcal(\tau,\tau^{\prime})}(a|\omega)\ln\paren{\frac{P_{\Tcal(\tau,\tau^{\prime})}(a|\omega)}{P_{\tau^{*}}(a|\omega_{\tau^{*}})}}\enspace .
\end{align*}
From \eqref{eq.dec.3} and \eqref{eq.dec.4}, for all $\epsilon>0$, we have therefore,
\begin{align*}
\sum_{\omega\in \Tcal(\tau,\tau^{\prime})}&\muh_{n-1}(\omega)\sum_{a\in A}(P_{\Tcal(\tau,\tau^{\prime})}(a|\omega)-\Ph(a|\omega))\ln\paren{\frac{P_{\tau^{\prime}}(a|\omega_{\tau^{\prime}})}{P_{\tau}(a|\omega_{\tau})}}\\
&\leq \epsilon\paren{K_{\mu}(P,\olP_{\tau})+K_{\mu}(P,\olP_{\tau^{\prime}})}\\
&\quad\quad\quad+(1+O(\eta))\frac{\Lst^{(\tau,\tau^{\prime})}}{\epsilon}\paren{K_{\mu_{\tau}}(P_{\tau},\Ph_{\tau})+K_{\mu_{\tau^{\prime}}}(P_{\tau^{\prime}},\Ph_{\tau^{\prime}})}\enspace .
\end{align*}
\end{proof}

\subsection{Upper bounds on $K_{\mu_{\tau}}$, $K_{\muh}$}
\begin{lemma}\label{lem.cons.typ}
Let $(X_{n})_{n\in\Z}$ be a stationary ergodic process satisfying assumption \eqref{hyp.conc}. Let $\delta>0$ and let $\Fcalst^{(n)}(\delta)$ and $\Fcalst(\delta)$ be the sets defined in \eqref{def.bon.ens.obs} and \eqref{def.bon.ens} respectively. Let $\Tbf_{typ}(\eta)$ be the set defined in \eqref{def.T.typ} and let $\Omega_{good}$ be the event \eqref{def.omega.good}. Then, on $\Omega_{good}$, for all $\tau\in \Fcalst^{(n)}(\delta)$, there exists $\eta=O\paren{\sqrt{\frac1{\Lambda_{n}^{(1)}}}\vee\frac1{\Lambda_{n}^{(2)}}}$ such that
\begin{align}
\label{eq.cont.K}K_{\mu_{\tau}}(P_{\tau},\Ph_{\tau})&\leq \paren{3+\eta}\paren{\sqrt{\rho_{n}}+\sqrt{ \frac{\varrho_{n}}{\Lambda_{n}^{(2)}}}}^{2}\sum_{(\omega,a)\in \tau\times A}\ln\paren{\frac1{\pi(\omega a)\delta}}\enspace .\\
\label{eq.cont.Kh}K_{\muh}(\Ph_{\tau},P_{\tau})&\leq\paren{3+\eta}\paren{\sqrt{\rho_{n}}+\sqrt{ \frac{\varrho_{n}}{\Lambda_{n}^{(2)}}}}^{2}\sum_{(\omega,a)\in \tau\times A}\ln\paren{\frac1{\pi(\omega a)\delta}}\enspace .
\end{align}
\end{lemma}

\begin{proof}
It comes from Lemma \ref{lem.cons.typ.gal} that, for some $\eta=O\paren{\sqrt{\frac1{\Lambda_{n}^{(1)}}}\vee\frac1{\Lambda_{n}^{(2)}}}$,
\begin{multline}\label{eq.cont.K.1}
K_{\mu}(\olP_{\tau},\Pt_{\tau})\leq\\
\paren{1+\eta}\sum_{(\omega,a)\in\tau\times A,\mu(\omega a)\neq 0}\frac{\paren{\muh_{n}(\omega a)-\mu(\omega a)}^{2}+2\paren{\muh_{n-1}(\omega)-\mu(\omega)}^{2}}{\mu(\omega)}.
\end{multline}
By Proposition \ref{prop.typ.gal}, $\Fcalst^{(n)}(\delta)\subset\Fcalst(\delta)$, hence, $\mu(\omega a)\geq \Lambda_{n}^{(2)}\varrho_{n}\ln(1/(\pi(\omega a)\delta))$. Thus
\begin{multline*}
\sqrt{\rho_{n}\mu(\omega a)\ln\paren{\frac1{\pi(\omega a)\delta}}}+\varrho_{n}\ln\paren{\frac1{\pi(\omega a)\delta}}\\\leq \paren{\sqrt{\rho_{n}}+\sqrt{ \frac{\varrho_{n}}{\Lambda_{n}^{(2)}}}}\sqrt{\mu(\omega a)\ln\paren{\frac1{\pi(\omega a)\delta}}}\enspace .
\end{multline*}
We obtain in the same way that
\begin{multline*}
\sqrt{\rho_{n}\mu(\omega)\ln\paren{\frac1{\pi(\omega)\delta}}}+\varrho_{n}\ln\paren{\frac1{\pi(\omega)\delta}}\\\leq \paren{\sqrt{\rho_{n}}+\sqrt{\frac{\varrho_{n}}{\Lambda_{n}^{(2)}}}}\sqrt{\mu(\omega)\ln\paren{\frac1{\pi(\omega)\delta}}}\enspace .
\end{multline*}
We deduce from assumption \eqref{hyp.conc} that
\begin{multline*}
\frac{\paren{\muh_{n}(\omega a)-\mu(\omega a)}^{2}+2\paren{\muh_{n-1}(\omega)-\mu(\omega)}^{2}}{\mu(\omega)}\\\leq (2+P(a|\omega))\paren{\sqrt{\rho_{n}}+\sqrt{ \frac{\varrho_{n}}{\Lambda_{n}^{(2)}}}}^{2}\ln\paren{\frac1{\pi(\omega a)\delta}}\enspace .
\end{multline*}
Plugging this last inequality in \eqref{eq.cont.K.1} yields \eqref{eq.cont.K}. \eqref{eq.cont.Kh} is obtained with the same inequality since, from Lemma \ref{lem.cons.typ.gal} we have for  $\eta=O\paren{\sqrt{\frac1{\Lambda_{n}^{(1)}}}\vee\frac1{\Lambda_{n}^{(2)}}}$,
\begin{multline*}
K_{\muh}(\Pt_{\tau},\olP_{\tau})\\\leq\paren{1+\eta}\sum_{(\omega,a)\in\tau\times A,\mu(\omega a)\neq 0}\frac{\paren{\muh_{n}(\omega a)-\mu(\omega a)}^{2}+2\paren{\muh_{n-1}(\omega)-\mu(\omega)}^{2}}{\mu(\omega)}.
\end{multline*}

\end{proof}

\subsection{Consequences of typicality}
\begin{lemma}\label{lem.cons.typ.gal}
Let $\eta<1/3$ and let $\tau\subset \Tbf_{typ}(\eta)$. We have
\begin{multline}
\label{eq.cons.typ.1}K_{\mu_{\tau}}(P_{\tau},\Ph_{\tau})\\\leq \paren{\frac12+\frac{2\eta}{3(1-3\eta)}}\sum_{(\omega,a)\in\tau\times A,\;\mu(\omega a)\neq 0}\mu(\omega)\frac{\paren{P_{\tau}(a|\omega)-\Ph(a|\omega)}^{2}}{P_{\tau}(a|\omega)}\enspace .
\end{multline}
Moreover, if we denote by 
\[
K_{\muh}(\Ph_{\tau},P_{\tau})\egaldef \sum_{(\omega,a)\in\tau\times A}\muh_{n}(\omega a)\ln\paren{\frac{\Ph_{\tau}(a|\omega)}{P_{\tau}(a|\omega)}}\enspace ,
\]
we have
\begin{multline}
\label{eq.cons.typ.2}K_{\muh}(\Ph_{\tau},P_{\tau})\\\leq \paren{\frac12+\frac{2\eta}{3(1-3\eta)}}\sum_{(\omega,a)\in\tau\times A,\;\muh_{n}(\omega a)\neq 0}\muh_{n-1}(\omega)\frac{\paren{P_{\tau}(a|\omega)-\Ph(a|\omega)}^{2}}{\Ph_{\tau}(a|\omega)}\enspace .
\end{multline}
In addition, if $\eta\rightarrow 0$,
\begin{equation}\label{eq.cons.typ.3}
\absj{K_{\mu_{\tau}}(P_{\tau},\Ph_{\tau})-K_{\muh}(\Ph_{\tau},P_{\tau})}=O(\eta)K_{\mu_{\tau}}(P_{\tau},\Ph_{\tau})\enspace .
\end{equation}
Finally, for all $(\omega,a)\in\tau\times A,\telque \mu(\omega a)\neq 0$,
\[
\mu(\omega)\frac{\paren{P_{\tau}(a|\omega)-\Ph(a|\omega)}^{2}}{P_{\tau}(a|\omega)}\leq 2\frac{(\muh_{n}(\omega a)-\mu(\omega a))^{2}}{\mu(\omega)}+4\frac{(\muh_{n-1}(\omega)-\mu(\omega))^{2}}{\mu(\omega)}\enspace .
\]
For all $(\omega,a)\in\tau\times A,\telque \muh_{n}(\omega a)\neq 0$,
\[
\muh_{n}(\omega)\frac{\paren{P_{\tau}(a|\omega)-\Ph(a|\omega)}^{2}}{\Ph_{\tau}(a|\omega)}\leq 2\frac{(\muh_{n}(\omega a)-\mu(\omega a))^{2}}{\mu(\omega)}+4\frac{(\muh_{n-1}(\omega)-\mu(\omega))^{2}}{\mu(\omega)}\enspace .
\]
\end{lemma}
\begin{proof}
Let us first remark that, for all $\omega$ in $\tau$ such that $\mu(\omega a)\neq 0$, we have,
\begin{equation}\label{eq.cons.typ}
\absj{\frac{\Ph_{\tau}(a|\omega)-P_{\tau}(a|\omega)}{P_{\tau}(a|\omega)}}\leq \frac{2\eta}{1-\eta}<1,\qquad \absj{\frac{\Ph_{\tau}(a|\omega)-P_{\tau}(a|\omega)}{\Ph_{\tau}(a|\omega)}}\leq \frac{2\eta}{1-\eta}<1\enspace .
\end{equation}
In addition, for all $\eta^{\prime}<1$, for all $u\leq \eta^{\prime}$, we have
\begin{equation}\label{eq.cont.ln}
\absj{-\ln(1-u)-u-\frac{u^{2}}2}\leq u^{2}\frac{\eta^{\prime}}{3(1-\eta^{\prime})}
\end{equation}
\eqref{eq.cons.typ.1}, \eqref{eq.cons.typ.2} and \eqref{eq.cons.typ.3} follow, plugging \eqref{eq.cons.typ} and \eqref{eq.cont.ln} in \eqref{eq.dec.Kull.1} and \eqref{eq.dec.Kull.2}, where, for $u=\frac{P_{\tau}(a|\omega)-\Ph(a|\omega)}{P_{\tau}(a|\omega)}$,
\begin{align}
\notag K_{\mu_{\tau}}&(P_{\tau},\Ph_{\tau})-\sum_{(\omega,a)\in\tau\times A,\;\mu(\omega a)\neq 0}\mu(\omega)\frac{(\Ph(a|\omega)-P_{\tau}(a|\omega))^{2}}{P_{\tau}(a|\omega)}\\
\notag&\egaldef \sum_{(\omega,a)\in\tau\times A}\mu(\omega a)\ln\paren{\frac{P_{\tau}(a|\omega)}{\Ph(a|\omega)}}\\
\notag&\quad\quad\quad-\sum_{(\omega,a)\in\tau\times A,\;\mu(\omega a)\neq 0}\mu(\omega)\frac{(\Ph(a|\omega)-P_{\tau}(a|\omega))^{2}}{P_{\tau}(a|\omega)}\\
\label{eq.dec.Kull.1}&=\sum_{(\omega,a)\in\tau\times A,\;\mu(\omega a)\neq 0}\mu(\omega a)\paren{-\ln\paren{1-u}-u-\paren{u}^{2}}\enspace .
\end{align}
And, for $u=\frac{\Ph_{\tau}(a|\omega)-P(a|\omega)}{\Ph_{\tau}(a|\omega)}$,
\begin{align}
\notag K_{\muh}&(\Ph_{\tau},P_{\tau})-\sum_{(\omega,a)\in\tau\times A,\;\muh_{n}(\omega a)\neq 0}\muh_{n}(\omega)\frac{(\Ph(a|\omega)-P_{\tau}(a|\omega))^{2}}{\Ph_{\tau}(a|\omega)}\\
\notag&\egaldef \sum_{(\omega,a)\in\tau\times A}\muh_{n}(\omega a)\ln\paren{\frac{\Ph_{\tau}(a|\omega)}{P_{\tau}(a|\omega)}}\\
\notag&\quad\quad\quad-\sum_{(\omega,a)\in\tau\times A,\;\muh_{n}(\omega a)\neq 0}\muh_{n}(\omega)\frac{(\Ph(a|\omega)-P_{\tau}(a|\omega))^{2}}{\Ph_{\tau}(a|\omega)}\\
\label{eq.dec.Kull.2}&=\sum_{(\omega,a)\in\tau\times A,\;\muh_{n}(\omega a)\neq 0}\muh_{n}(\omega a)\paren{-\ln\paren{1-u}-u-\paren{u}^{2}}\enspace .
\end{align}
The bounds on \(
\mu(\omega)\frac{\paren{P_{\tau}(a|\omega)-\Ph(a|\omega)}^{2}}{P_{\tau}(a|\omega)}
\) follow from the inequalities
\[
\absj{P_{\tau}(a|\omega)-\Ph(a|\omega)}\leq\frac{\absj{\mu(\omega a)-\muh_{n}(\omega a)}}{\mu(\omega)}+\Ph(a|\omega)\frac{\absj{\mu(\omega)-\muh_{n-1}(\omega)}}{\mu(\omega)}\enspace .
\]
\[
\absj{P_{\tau}(a|\omega)-\Ph(a|\omega)}\leq\frac{\absj{\mu(\omega a)-\muh_{n}(\omega a)}}{\muh_{n-1}(\omega)}+P(a|\omega)\frac{\absj{\mu(\omega)-\muh_{n-1}(\omega)}}{\muh_{n-1}(\omega)}\enspace .
\]
These imply in particular, since $\eta\leq 1/3$,
\begin{align*}
\absj{P_{\tau}(a|\omega)-\Ph(a|\omega)}&\leq\frac{\absj{\mu(\omega a)-\muh_{n}(\omega a)}}{\sqrt{\mu(\omega)\muh_{n-1}(\omega)}}+\\
\notag&\quad\quad\quad\paren{\Ph(a|\omega)\vee P(a|\omega)}\frac{\absj{\mu(\omega)-\muh_{n-1}(\omega)}}{\sqrt{\mu(\omega)\muh_{n-1}(\omega)}}\\
&\leq \frac{\absj{\mu(\omega a)-\muh_{n}(\omega a)}}{\sqrt{\mu(\omega)\muh_{n-1}(\omega)}}+\sqrt{2P(a|\omega)}\frac{\absj{\mu(\omega)-\muh_{n-1}(\omega)}}{\sqrt{\mu(\omega)\muh_{n-1}(\omega)}}\enspace .
\end{align*}
\end{proof}

\subsection{Tools for mixing processes}

\begin{lemma}\label{lem.covariance}
Let $(X_{n})_{n\in \Z}$ be a $\phi$-mixing process satisfying \eqref{cond.phi.mix} and \eqref{cond.ND}. Then, for all $\omega\in A^{*}$,
\begin{equation}\label{eq.cov.1}
\sum_{k=0}^{\infty}\absj{\cov\paren{\un_{X_{-|\omega|+1}^{0}=\omega},\un_{X_{k-\absj{\omega}+1}^{k}=\omega}}}\leq \paren{\Phi+\frac1{1-\lambda}}\mu(\omega)\enspace .
\end{equation}
As a consequence, for all $N\in \N^{*}$, 
\[
\var\paren{\sum_{k=0}^{N}\un_{X_{k-\absj{\omega}+1}^{k}=\omega}}\leq 2N\paren{\Phi+\frac1{\lambda-1}}\mu(\omega)\enspace .
\]
If, in addition, for any $(a,\omega)\in A\times A^{*}$, $P(X_{0}=a|X_{1}^{|\omega|}=\omega)\leq \lambda$, then, for all $(\omega,\omega^{\prime})\in( A^{*})^{2}$,
\begin{multline*}
\absj{\cov\paren{\un_{X_{-|\omega^{\prime}|+1}^{0}=\omega^{\prime}},\un_{X_{k-\absj{\omega}+1}^{k}=\omega}}}\\\leq \paren{\sqrt{\phi_{k-\absj{\omega}+1}}\un_{k-\absj{\omega}+1\geq 0}+\lambda^{k}\un_{k-\absj{\omega}+1< 0}}\sqrt{\mu(\omega)\mu(\omega^{\prime})}\enspace .
\end{multline*}
\end{lemma}
\begin{proof}
If $k-\absj{\omega}+1\geq 0$, we use Lemma \ref{lem.cov.Vi97} and we get
\[
\absj{\cov\paren{\un_{X_{-|\omega^{\prime}|+1}^{0}=\omega^{\prime}},\un_{X_{k-\absj{\omega}+1}^{k}=\omega}}}\leq \sqrt{\phi_{k-\absj{\omega}+1}\mu(\omega)\mu(\omega^{\prime})}\enspace .
\]
Therefore,
\begin{multline}\label{eq.cov.int.1}
\sum_{k=\absj{\omega}-1}^{\infty}\absj{\cov\paren{\un_{X_{-|\omega^{\prime}|+1}^{0}=\omega^{\prime}},\un_{X_{k-\absj{\omega}+1}^{k}=\omega}}}\\\leq \sum_{k=0}^{\infty}\sqrt{\phi_{k}\mu(\omega)\mu(\omega^{\prime})}\leq \Phi\sqrt{\mu(\omega)\mu(\omega^{\prime})}\enspace .
\end{multline}
If $k-\absj{\omega}+1< 0$, denoting by $\omega=a_{-\absj{\omega}}^{-1}$, condition \eqref{cond.ND} implies
\begin{align*}
&\absj{\cov\paren{\un_{X_{-|\omega^{\prime}|+1}^{0}=\omega^{\prime}},\un_{X_{k-\absj{\omega}+1}^{k}=\omega}}}\leq\mu(\omega^{\prime})\absj{\mu\paren{X_{1}^{k}=a^{-1}_{-k}|X_{-\absj{\omega^{\prime}}}^{0}=\omega^{\prime}}}\\
&\leq \mu(\omega^{\prime})\prod_{l=1}^{k}\absj{\mu\paren{X_{k}=a_{-1}|X_{-\absj{\omega^{\prime}}}^{k-1}=\omega^{\prime} a_{-k}^{-2}}}\leq \mu(\omega^{\prime})\lambda^{k}\enspace .
\end{align*}
This is sufficient to obtain \eqref{eq.cov.1}, choosing $\omega^{\prime}=\omega$. As a consequence,
\begin{align*}
\var&\paren{\sum_{k=0}^{N}\un_{X_{k-\absj{\omega}+1}^{k}=\omega}}=\sum_{k,k^{\prime}=0}^{N}\cov\paren{\un_{X_{k-\absj{\omega}+1}^{k}=\omega},\un_{X_{k^{\prime}-\absj{\omega}+1}^{k^{\prime}}=\omega}}\\
&\leq 2\sum_{k=0}^{N}(N-k+1)\absj{\cov\paren{\un_{X_{-|\omega|+1}^{0}=\omega},\un_{X_{k-\absj{\omega}+1}^{k}=\omega}}}\\
&\leq 2N\sum_{k=0}^{\infty}\absj{\cov\paren{\un_{X_{-|\omega|+1}^{0}=\omega},\un_{X_{k-\absj{\omega}+1}^{k}=\omega}}}\leq 2N\paren{\Phi+\frac1{\lambda-1}}\mu(\omega)\enspace .
\end{align*}
An argument symmetric to the one in \eqref{eq.cov.int.1} shows that, if $k-\absj{\omega}+1< 0$, when we have moreover, for any $(a,\omega)\in A\times A^{*}$, $P(X_{0}=a|X_{1}^{|\omega|}=\omega)\leq \lambda$, 
\[
\absj{\cov\paren{\un_{X_{-|\omega^{\prime}|+1}^{0}=\omega^{\prime}},\un_{X_{k-\absj{\omega}+1}^{k}=\omega}}}\leq \lambda^{k}\mu(\omega)\enspace .
\]
Thus, 
\[
\absj{\cov\paren{\un_{X_{-|\omega^{\prime}|+1}^{0}=\omega^{\prime}},\un_{X_{k-\absj{\omega}+1}^{k}=\omega}}}\leq \lambda^{k}\mu(\omega)\wedge \mu(\omega^{\prime})\leq \lambda^{k}\sqrt{\mu(\omega)\mu(\omega^{\prime})}\enspace .
\]
\end{proof}

The following lemmas are due to Viennet \cite{Vi97}

\begin{lemma}\label{lem.Vi97}
Let $(X_{n})_{n\in \Z}$ be a $\beta$-mixing process. Let $(J_{k})_{k=1,\ldots,N}$ be a collection of subsets of $\N$ satisfying the following conditions.
\begin{enumerate}
\item $\exists \;q_{o}\in \N^{*}$ such that, for all $k=1,\ldots,N-1$, $\max\set{i\in J_{k}}\leq q_{o}+\min\set{j\in J_{k+1}}.$
\item $\exists \;M_{o}\in \N^{*}$ such that, for all $k=1,\ldots,N$, $\card\set{J_{k}}\leq M_{o}$.
\end{enumerate}
Then, there exists random variables $(Y_{i})_{i\in \cup_{k=1}^{N}J_{k}}$ such that,
\begin{enumerate}
\item for all $k=1,\ldots,N$, $(Y_{i})_{i\in J_{k}}$ has the same distribution as $(X_{i})_{i\in J_{k}}$,
\item for all $k=2,\ldots,N$, $(Y_{i})_{i\in J_{k}}$ is independent of $(X_{i},Y_{i})_{i\in \cup_{t\leq k-1}J_{t}}$,
\item for all $k=1,\ldots,N$, $\P\set{(X_{i})_{i\in J_{k}}\neq (Y_{i})_{i\in J_{k}}}\leq M_{o}\beta_{q_{o}}$.
\end{enumerate}
\end{lemma}

\begin{lemma}\label{lem.cov.Vi97}
Let $X$, $Y$ be two real valued random variables. There exists two real functions $b_{1}$ and $b_{2}$ such that, for all bounded functions $f$ and $g$,
\begin{align}
\label{cond.mix} \norm{b_{1}}_{\infty}&\leq \phi(\sigma(X),\sigma(Y)),\;\norm{b_{2}}_{\infty}\leq \phi(\sigma(Y),\sigma(X))\enspace .\\
\label{cov.ineq} \cov\paren{f(X),g(Y)}\leq &\sqrt{\E\paren{b_{1}(X)f^{2}(X)}\E\paren{b_{2}(Y)g^{2}(Y)}}\enspace .
\end{align}
\end{lemma}

\subsection{Additional lemmas}
The following lemma can be found, for example, in \cite{Ma07} Lemma 7.24.
\begin{lemma}\label{lem:BarSheu}
For all probability measures $P$, $Q$ with $P<<Q$,
\[
\frac12\int(dP\wedge dQ)\paren{\ln\paren{\frac{dP}{dQ}}}^{2}\leq \int dP\ln\paren{\frac{dP}{dQ}}\enspace .
\]
\end{lemma}

\begin{lemma}\label{lem.pyth.gal}
Let $\mu$ be a probability measure with kernel $P$ and let $\tau$ be a finite tree. Then, for all $\nu\in \Mcal_{\tau}$ with transition kernel $Q$ such that $K_{\mu}(P,Q)<\infty$, we have
\[
K_{\mu}(P,Q)=K_{\mu}(P,\olP_{\tau})+K_{\mu_{\tau}}(P_{\tau},Q_{\tau})\enspace .
\]
\end{lemma}
\begin{proof}
By definition,
\begin{align*}
K_{\mu}\paren{P,Q}&=\int_{A^{-\N}\times A}d\mu(\omega a)\ln\paren{\frac{P(a|\omega)}{Q(a|\omega)}}\\
=&\int_{A^{-\N}\times A}d\mu(\omega a)\ln\paren{\frac{P(a|\omega)}{\olP_{\tau}(a|\omega)}}+\int_{A^{-\N}\times A}d\mu(\omega a)\ln\paren{\frac{\olP_{\tau}(a|\omega)}{Q(a|\omega)}}\\
=&K_{\mu}\paren{P,\olP_{\tau}}+\int_{A^{-\N}\times A}d\mu(\omega a)\ln\paren{\frac{\olP_{\tau}(a|\omega)}{Q(a|\omega)}}\enspace .
\end{align*}

For all $\omega\in A^{-\N}$, let $\omega_{1}\in A^{-\N}$ such that $\omega=\omega_{1}\omega_{\tau}$. As the function $\ln\paren{\frac{\olP_{\tau}(a|\omega)}{Q(a|\omega)}}$ does not depend on $\omega_{1}$, is equal to $\ln\paren{\frac{P_{\tau}(a|\omega_{\tau})}{Q_{\tau}(a|\omega_{\tau})}}$ and $\mu$ satisfies, for all $(\omega,a)\in \tau\times A$, $\int_{\omega_{1}\in A^{-\N}}d\mu(\omega_{1}\omega a)=\mu(\omega a)=\mu_{\tau}(\omega a)$, we have

\begin{multline*}
\int_{A^{-\N}\times A}d\mu(\omega a)\ln\paren{\frac{\olP_{\tau}(a|\omega)}{Q(a|\omega)}}=\sum_{(\omega_{\tau},a)\in \tau\times A}\mu_{\tau}(\omega_{\tau} a)\ln\paren{\frac{P_{\tau}(a|\omega_{\tau})}{Q_{\tau}(a|\omega_{\tau})}}\\=K_{\mu_{\tau}}\paren{P_{\tau}, Q_{\tau}}\enspace .
\end{multline*}
\end{proof}

\begin{lemma}\label{lem.benett} (Benett's inequality)
Let $\xi_{1:N}$ be independent random variables such that, $\forall i=1,\ldots,N$, $\norm{\xi_{i}}_{\infty}\leq b$. Then, for all $y>0$, 
\[
\P\set{\sum_{i=1}^{N}(Y_{i}-\E(Y_{i})\geq \sqrt{2\sum_{i=1}^{N}\var(Y_{i})y}+\frac{by}3}\leq e^{-y}\enspace .
\]
\end{lemma}

% AOS,AOAS: If there are supplements please fill:
%\begin{supplement}[id=suppA]
%  \sname{Supplement A}
%  \stitle{Title}
%  \slink[url]{http://lib.stat.cmu.edu/aoas/???/???}
%  \sdescription{Some text}
%\end{supplement}

\end{document}